\documentclass[reqno]{amsart}
\usepackage{latexsym}
\usepackage{amstext}
\usepackage {amsmath}
\usepackage {amsfonts}
\usepackage {amssymb}
\usepackage {amsthm}
\usepackage {bbm}
\usepackage{enumerate}
\usepackage{array}
\usepackage{comment}
\usepackage{xspace}
\usepackage{graphicx}
\def\loc{{\mathrm{loc}}}

\DeclareMathOperator{\lin}{lin}

\newcommand{\eps}{\varepsilon}

\newcommand{\R}{\ensuremath{\mathbb{R}}}

\newcommand{\N}{\ensuremath{\mathbb{N}}\ }

\newcommand{\K}{\ensuremath{\mathcal K}}

\newcommand{\schwto}{\ensuremath{\to}}

\newcounter{counter-liste}
\newcounter{counter-liste2}

  {\ \begin{list}{{(\arabic{counter-liste})}\hfill}%
  {\topsep2mm\itemindent1ex\leftmargin0cm\usecounter{counter-liste}}
  }%
  {\end{list}}
\newenvironment{liste(a)}%
  {\ \begin{list}{{(\alph{counter-liste})}\hfill}%
  {\topsep2mm\itemindent1ex\leftmargin0cm\usecounter{counter-liste}}
  }%
  {\end{list}}
\def\R{\mathbb R}

\theoremstyle{plain}
\numberwithin{equation}{section}
\newtheorem{lemma}{Lemma}[section]
\newtheorem{theorem}[lemma]{Theorem}
\newtheorem{proposition}[lemma]{Proposition}
\newtheorem{definition}[lemma]{Definition}

\newtheorem{corollary}[lemma]{Corollary}

\theoremstyle{definition}
\newtheorem{remark}[lemma]{Remark}
\begin{document}

\title[Travelling waves in a tissue-degradation model]{Travelling-wave analysis
  of a model describing tissue degradation by bacteria} 
\thanks{
The research of the authors was supported by
the European Community's
Human Potential Programme
under contract HPRN-CT-2002-00274,
FRONTS-SINGULARITIES.
}

\author[Hilhorst]{Danielle Hilhorst}
\address{Danielle Hilhorst, CNRS and Laboratoire de Math{\'e}matiques,
  Universit{\'e} de Paris-Sud, 91405 Orsay Cedex, France}

\author[King]{John R. King}
\address{John R. King, Centre for Mathematical Medicine, Theoretical Mechanics Section,
School of Mathematical \mbox{Sciences},
University of Nottingham,
Nottingham, NG7 2RD, UK}

\author[R\"{o}ger]{Matthias R\"{o}ger}
\address{Matthias R\"{o}ger, Max Planck Institute for Mathematics in the Sciences,
Inselstr. 22,
D-04103 Leipzig}



\date{\today}

\begin{abstract}
We study travelling-wave solutions for a reaction-diffusion
system arising as a model for host-tissue degradation by bacteria.
This system consists of a parabolic equation coupled with an
ordinary differential equation.
For large values of the `degradation-rate parameter' solutions are
well approximated by solutions of a Stefan-like free boundary
problem, for which travelling-wave solutions can be found explicitly.
Our aim is to prove the existence of travelling waves
for all sufficiently large wave-speeds for the original reaction-diffusion system
and to determine the minimal speed.
We prove that for all sufficiently large degradation rates the minimal speed
is identical to
the minimal speed of the limit problem. In particular, in this parameter range,
{\itshape nonlinear} selection of the minimal speed occurs.
\end{abstract}
\maketitle
\tableofcontents
%=============================================
% introduction
%=============================================
\section{Introduction}
In this article we continue our mathematical analysis of a model
for the degradation of host tissue by extracellular bacteria.
This model was introduced in \cite{Ki} and consists of a
reaction-diffusion equation coupled with an ordinary differential equation.
In \cite{HKR} we proved the existence
of solutions to the time-dependent problem and the convergence to a
limit problem in the  
`large-degradation-rate' limit.
Here we turn to the question of existence
and behaviour of travelling-wave solutions.

There is an increasing interest in models which support the understanding of bacterial
infections, and we refer to \cite{Ki} and \cite{Wa} for further background and references.
This paper is in effect concerned with the specific issue of how rapidly
a bacterial infection in, for example, a burn wound may invade
the underlying
tissue (with dire potential consequences for the patient, notably mortality due
to septicemia). For the type of model with which we are concerned
here, the relevant invasion speed is expected to be governed by the corresponding
travelling-wave problem. Typically, the smallest possible wave-speed is
realized by a large class of solutions. Accordingly determining the
minimal speed of travelling waves becomes a central question
(with obvious implications
for the amount of time available for medical treatment, for instance).

In a dimensionless form, the model in \cite{Ki} is given by the equations
\begin{eqnarray}
        \partial_t {u} &=& \Delta{u} - {u} + {w} -\gamma k {u}(1-{w}),
                \label{rd-u}\\
        \partial_t {w} &=& k {u}(1-{w}),\label{rd-w}
\end{eqnarray}
where $u$ describes the concentration of degradative enzymes, $(1-w)$
the volume fraction of healthy tissue and $\gamma,k$ are positive constants.
The key parameter here is the \emph{ degradation-rate} $k$, which is  very large
in practice.
Equations \eqref{rd-u}, \eqref{rd-w} are considered in a time-space cylinder,
with the upper half space of $\R^3$ as the spatial domain.
Finally, the system is complemented by initial conditions for $u$ and $w$,
a Neumann condition on the lateral boundary for $u$ and a decay condition
for $u$ and $w$ in the far field. In \cite{HKR} we gave a precise mathematical
formulation and proved the existence and uniqueness of solutions to a slightly more general
system, including the possibility of a diffusion term in \eqref{rd-w}.
One noteworthy aspect of \eqref{rd-u}, \eqref{rd-w} is the convergence of solutions
to the solution of a Stefan-like free boundary problem as the degradation rate $k$
tends to infinity. This \emph{large-degradation-rate limit} was identified
by a formal asymptotic analysis in \cite{Ki} and was proved in \cite{HKR}.

Reaction-diffusion systems of the general form
\begin{eqnarray}\label{eq-rd-gen}
        V_t &=& D\Delta V + F(V),
\end{eqnarray}
where $V:(0,T)\times\R^n\to \R^N$ is vector-valued, $F:\R^N\to\R^N$ is a given nonlinearity
and $D\in\R^{N\times N}$ is a diagonal positive-semi-definite matrix,
appear in a lot of different scientific areas.
One-dimensional travelling waves are solutions on $(0,\infty)\times\R$ of the special form
\begin{eqnarray*}
  \tilde{V}(t,x) &=& V(x-ct),
\end{eqnarray*}
where $c$ is called the speed and $V$ the profile of this travelling wave.
The question of existence and behaviour of travelling waves is of
enormous interest in many of the applications and pertinent results for
the vector case $N\geq 2$ remain restricted to rather specific systems.

The system \eqref{rd-u}, \eqref{rd-w} has, as we will see in Remark
\ref{rem:mono}, one stable equilibrium in $(u,w)=(1,1)$ and one
unstable equilibrium in $(u,w)=(0,0)$ and therefore belongs to the class of
monostable systems. Scalar monostable equations
\begin{eqnarray}\label{eq-Fi}
  v_t &=& \Delta v + f(v)
\end{eqnarray}
where
\begin{gather*}
  f(0)\,=\,f(1)\,=\,0,\quad f(v)>0\text{ for }0<v<1,\quad f^\prime(0)>0>f^\prime(1),
\end{gather*}
are well-studied, especially the famous Fisher equation, that is \eqref{eq-Fi} with
$f(v)=v(1-v)$, introduced in \cite{Fi}.
The rigorous analysis of equations of this type also started in the 1930s with the work of
Kolmogorov, Petrovskii and Piskunov \cite{KPP}. Under an extra assumption on $f$ they
proved the existence of travelling waves for all speeds $c\geq c_{\lin}$, where
$c_{\lin}$ can be found explicitly in terms of $f$ by a linearisation about $v=0$
(corresponding to a degenerate node in the travelling-wave phase plane).
Moreover, they proved that the solutions to \eqref{eq-Fi} with initial data decaying sufficiently fast
propagate with speed $c_{\lin}$.
For more general monostable $f$ the propagation speed was found to be either equal to
or larger than $c_{\lin}$ and therefore one distinguishes between a {\itshape linear}
or {\itshape nonlinear} selection of the propagation speed (the terminology
{\itshape pulled} and {\itshape pushed} fronts, respectively, having an equivalent
meaning).
Aronson and Weinberger \cite{AW} (see also Hadeler and Rothe \cite{HR} and \cite{Sto}
for other pioneering work on such matters and \cite{Saa-rev} for a recent review)
proved that for general monostable $f$, in both
the linear and the nonlinear selection case, the propagation speed for
solutions with initial data decaying sufficiently fast is given by
the minimal speed of travelling waves. They showed that monotonic travelling waves exist
for all speeds $c\geq c_{\min}$ and none for $c<c_{\min}$, where $c_{\min}\geq c_{\lin}$;
solutions of \eqref{eq-Fi} with sufficiently rapidly
decaying initial data propagate with speed $c_{\min}$.
For the nonlinear selection cases, $c_{\min}>c_{\lin}$, Rothe \cite{Rot} and Roquejoffre \cite{Roq}
proved that, if the initial data decays sufficiently rapidly, the large-time solutions to
\eqref{eq-Fi} not
only propagate with speed $c_{\min}$ but also approach the profile of a travelling wave
with minimal speed.

Whereas the connection between large-time behaviour and existence of travelling waves
for monostable equations is satisfactorily resolved, the distinction between nonlinear
or linear selection is still a challenging question, see for example
\cite{Saa}, \cite{ES}, \cite{BD}.
Only a few rigorous results for general monostable nonlinearities are available.
In \cite{LMN} a variational characterisation of travelling waves and
a concrete criterion for whether linear or nonlinear selection occurs for a given equation
was derived.
Even fewer general analytical results are available for monostable {\itshape systems}. An
existence theorem for travelling waves was proved in \cite{Vo} for monotone monostable systems,
which are systems of the form  \eqref{eq-rd-gen} in which the Jacobian matrix
$DF$ has only nonnegative off-diagonal elements.
Results on the existence of travelling waves and the long time behaviour of \eqref{eq-rd-gen}
for monostable gradient systems, that is for $n=N$ and nonlinearities $F=\nabla g$ with
$g:\R^n\to\R$, were given in \cite{Mur}.
To the best of our knowledge there are no more general results on the question of whether
linear or nonlinear selection will occur.

In this article we prove that for all $\gamma,k$ there exist
monotone travelling waves for the system \eqref{rd-u}, \eqref{rd-w}
for all speeds $c\geq c_{\min}$ and no speeds $c< c_{\min}$.
The minimal speed $c_{\min}$ in general depends on the
parameters $\gamma,k$.

We prove that for all $k>k_0$, where $k_0$ is explicitly given
in terms of $\gamma$, the minimal speed is larger than the value $c_{\lin}(k)$ obtained from
a linearisation at the unstable equilibrium. Surprisingly enough, for $k>k_0$ the minimal
speed of travelling waves is identical to the minimal speed of travelling waves for the
Stefan-like limit problem that was formulated in \cite{HKR}.
Our analysis is based on two main facts. One is the monotone structure of the
system, which makes possible the use of comparison principles for the parabolic problem.
The second is a remarkable reduction in order of the travelling-wave equations
that occurs when $c$ is given by the minimal speed $c_\infty$ of travelling waves of
the large-degradation limit.
We obtain the existence of travelling waves with speed $c_\infty$ for \eqref{rd-u}, \eqref{rd-w}.
Finally, a comparison argument allows us to prove that nonlinear selection occurs
for sufficiently large values of $k$, the minimal speed in this regime being identical to $c_\infty$.

This paper is organised as follows. In section \ref{sec-exist} we prove the existence of travelling
waves for the reaction-diffusion system \eqref{rd-u}, \eqref{rd-w}. In section \ref{sec-limit}
we recall the formulation of the large degradation limit problem
and consider travelling-wave solutions for this problem.
In section \ref{sec-speed} we return to the reaction-diffusion system and investigate
the selection of the minimal speed. Section \ref{sec-conv} deals with the convergence
of travelling waves for the reaction-diffusion system to travelling waves of the Stefan-like
free boundary problem as the reaction rate $k$ approaches infinity.
Finally, we give some conclusions and remarks on open problems in section \ref{sec-concl}.
%=====================================================
%       section existence
%=====================================================
\section{Existence of monotone travelling waves}
\label{sec-exist}
In this section we prove the existence
of monotone travelling waves. First we fix some notation and make some
remarks.
\begin{remark}\label{rem:mono}
A system of the general form
\eqref{eq-rd-gen} is called \emph{monotone} if the off-diagonal elements
of the Jacobian matrix $DF$ are
non-negative and \emph{strictly monotone} if they are positive (see \cite{Vo}).
A system of the form \eqref{eq-rd-gen} with two stationary points
is called \emph{monostable} if one of the stationary points is stable
and the other is unstable.\\
The system \eqref{rd-u}, \eqref{rd-w} is of the form \eqref{eq-rd-gen}
with
\begin{gather*}
  V\,=\,
  \begin{pmatrix}
    u\\
    w
  \end{pmatrix}\!,\qquad
  F(z_1,z_2)\,=\, 
  \begin{pmatrix}
    -z_1 +z_2 -\gamma kz_1(1-z_2)\\
    kz_1(1-z_2)
  \end{pmatrix}\!.
\end{gather*}
It follows that \eqref{rd-u}, \eqref{rd-w} is a monotone but not strictly monotone system.
Further, $z=(0,0)$ and $z=(1,1)$ are the only stationary points of \eqref{rd-u}, \eqref{rd-w}
and we obtain
\begin{gather*}
  DF(0,0)\,=\,
  \begin{pmatrix}
    -(1+\gamma k) & 1\\
    k & 0
  \end{pmatrix}
\end{gather*}
with one positive and one negative eigenvalue, and
\begin{gather*}
  DF(1,1)\,=\, 
  \begin{pmatrix}
    -1 & 1+\gamma k\\
    0 & -k
  \end{pmatrix}
\end{gather*}
with two negative eigenvalues. Therefore \eqref{rd-u}, \eqref{rd-w}
is a monostable monotone system.
\end{remark}
As remarked before, for a one-dimensional travelling wave $u,w$ of \eqref{rd-u},
\eqref{rd-w} with speed $c$, the functions
\begin{gather*}
        \tilde{u}(t,x) \,:=\, u(x-ct),\quad
        \tilde{w}(t,x) \,:=\, w(x-ct),
\end{gather*}
are solutions of \eqref{rd-u}, \eqref{rd-w} on $(0,\infty)\times\R$.
Therefore $u,w$ have to satisfy the \emph{travelling-wave equations}
\begin{eqnarray}
        0 &=& u^{\prime\prime} + cu^\prime -u + w -\gamma k u(1-w),\label{eq-tw-u}\\
        0 &=& cw^\prime +k u(1-w),\label{eq-tw-w}
\end{eqnarray}
We restrict our investigations to functions $u,w$ taking values only in $[0,1]$,
which is the range of
meaningful values in the tissue degradation model. % Observe that $(u,w)=(0,0)$ and $(u,w)=(1,1)$
%are the only stationary solutions of \eqref{eq-tw-u}, \eqref{eq-tw-w}.
We first summarise some basic properties of travelling-wave solutions.
\begin{lemma}\label{lem-tw-basic}
Assume that $u\in C^2(\R), w\in C^1(\R)$ with $0\leq u,w\leq 1$ satisfy
\eqref{eq-tw-u}, \eqref{eq-tw-w} with $c\in\R$.
Then the following properties hold.
\begin{enumerate}[(1)\hspace{-3pt}]
\item Invariance under space shifts. For any $x_0\in\R$ the functions defined by
$x\mapsto u(x+x_0)$, $x\mapsto w(x+x_0)$
solve \eqref{eq-tw-u}, \eqref{eq-tw-w} with the same $c\in\R$.
\item Invariance under inversion. The functions defined by $x\mapsto u(-x)$, $x\mapsto w(-x)$
solve \eqref{eq-tw-u}, \eqref{eq-tw-w} with $c$ replaced by $-c$.
\item Smoothness of travelling waves. The functions $u,w$ are infinitely differentiable.
\item Monotonicity of travelling waves. If $c>0$ and $u,w$ are not both constant
then $u,w$ are strictly monotone decreasing
and they approach unity as $x\to -\infty$ and zero
as $x\to\infty$.
\item Standing waves. For $c=0$ either $(u,w)=(0,0)$ or $(u,w)=(1,1)$.
\end{enumerate}
\end{lemma}
\begin{proof}
The statements (1), (2) are immediate from \eqref{eq-tw-u},
\eqref{eq-tw-w}. Statement (3) follows from a bootstrapping argument.\\
To prove (4) we observe from \eqref{eq-tw-w} that
$w^\prime\leq 0$ and $w^\prime <0$ if $u(1-w)>0$.
Next we differentiate \eqref{eq-tw-u} and obtain
\begin{eqnarray}
 0 &=& u^{\prime\prime\prime} + cu^{\prime\prime} -\big(1+\gamma k (1- w)\big)u^\prime
 +(1+\gamma ku) w^\prime.\label{eq:diff-u-eq}
\end{eqnarray}
From this equation we see that $u^\prime$ cannot have a positive local maximum.
Since $u$ is bounded,
$u^\prime$ cannot approach a positive supremum at $+\infty$ or $-\infty$
and we obtain that $u^\prime\leq 0$.

Assume now that $u^\prime(x_0)=0$. Therefore $u^\prime(x_0)$ is a local
maximum of $u^\prime$ and we deduce from \eqref{eq:diff-u-eq} that
$w^\prime(x_0)\geq 0$ and that, using $w^\prime\leq 0$ and
\eqref{eq-tw-w},
\begin{gather*}
  -\frac{k}{c} u(1-w)(x_0) \,=\, w^\prime(x_0)\,=\,0.
\end{gather*}
Let us consider the case  $u(x_0)=0$
(the case that $w$ takes the value one is analogous). Then, by $u^\prime\leq 0$ and
$u\geq 0$, we obtain
$u=0$ on $[x_0,\infty)$. From \eqref{eq-tw-u} it also follows that $w=0$ on $[x_0,\infty)$.
Since $u,w$ are smooth we find that $u,w$ solve \eqref{eq-tw-u}, \eqref{eq-tw-w}
with $u(x_0)=u^\prime(x_0)=w(x_0)=0$.
On the other hand solutions of \eqref{eq-tw-u}, \eqref{eq-tw-w} with these data prescribed
at $x_0$ are unique and thus identically equal to zero, in
contradiction to our assumption of non-constant solutions. Therefore $u(1-w)>0$ and $u^\prime,w^\prime<0$
holds.\\
Next, since $u(x),w(x)$ are monotone in $x$ and uniformly bounded, their limits as $x\to -\infty$ exist.
From \eqref{eq-tw-w} we deduce that $w^\prime(x)$ has a limit and, from $0\leq w\leq 1$, that
%\mbox{$w^\prime(x_j)\to 0$} for a subsequence $x_j\to-\infty$, which ensures
\begin{gather}\label{limit-u(1-w)}
        w^\prime(x) \,\to\, 0,\quad
        u(x)(1-w(x))\,\to\, 0\quad \text{ as }x\to-\infty.
\end{gather}
From \eqref{eq-tw-u} we similarly obtain that
\begin{eqnarray*}
        (u^\prime+cu)^\prime(x) &\to& 0\quad \text{ as }x\to-\infty,
\end{eqnarray*}
and thus
\begin{eqnarray*}
        -u(x)+w(x) &\to& 0 \quad\text{ as }x\to-\infty.
\end{eqnarray*}
By $0\leq u,w\leq 1$, \eqref{limit-u(1-w)}, and $u^\prime, w^\prime\leq 0$ this gives
\begin{gather*}
        u(x) \,\to\,1,\quad w(x) \,\to\, 1 \qquad\text{ as }x\to-\infty.
\end{gather*}
The proof that $u(x),w(x)$ approach zero as $x\to +\infty$ is similar.
%We now show that $u^\prime$ tends to zero as $x\to -\infty$. Observe that this function is bounded due to
%\eqref{bounds-u,w_x-infty} and assume the existence of two different limit points.
%Then we find a value $d<0$ in between and
%a subsequence $x_j\to -\infty $ such that $u^\prime(x_j)=d, u^{\prime\prime}(x_j)\leq 0$,
%in contradiction to \eqref{limit-u-prime+}. Therefore
%\begin{eqnarray*}
%       u^\prime(x) &\to& 0\quad\text{ as }x\to-\infty.
%\end{eqnarray*}

To prove claim (5) we first deduce from $c=0$ and \eqref{eq-tw-w} that
\begin{gather}
  u(1-w)\,=\, 0. \label{eq:u(1-w)}
\end{gather}
Let us assume that $u(x_0)>0$ for a point $x_0\in \R$. Then there exists
a maximal interval 
$(a,b)$ such that $x_0\in (a,b)$ and $u>0$ in $(a,b)$. From
\eqref{eq-tw-u}, \eqref{eq:u(1-w)}, and $c=0$ we then  obtain that 
\begin{gather}
  u\,>\, 0,\, w\,=\,1,\, u^{\prime\prime}\,<\, 0\qquad\text{ in }(a,b).
  \label{eq:prop-ab} 
\end{gather}
If $a>-\infty$ we deduce from the regularity assumptions on $u,w$ and
$u\geq 0$ that
\begin{gather*}
  u(a)\,=\,0,\, w(a)\,=\,1,\, u^\prime(a)=0
\end{gather*}
and finally $u^{\prime\prime}(a)\geq 0$, which yields a contradiction to
\eqref{eq-tw-u}. This shows that $a=-\infty$; by analogous arguments we
obtain that $b=\infty$ and $(a,b)=\R$. Therefore
\begin{gather*}
  0\,=\,u^{\prime\prime} -u+1\quad\text{ on }\R.
\end{gather*}
and from the boundedness of $u$ it follows that $u=1=w$.

Similarly we prove that $w=0=u$ if there exists $x_0\in\R$ such that
$w(x_0)<1$.
\end{proof}
Due to Lemma \ref{lem-tw-basic} we can restrict our investigations to
the following set of {\itshape admissible functions}.
\begin{definition}\label{def-admiss}
Let $\K$ be the set of functions $v\in C^\infty(\R)$, with
\begin{gather}
        0\,<\,v\,<\,1,\quad
        v^\prime\,<\,0,\label{K-1}\\
        v(x)\,\to\, 1,\quad \text{ as } x\to -\infty,\label{K-2}\\
        v(x)\,\to\, 0,\quad \text{ as } x\to +\infty.\label{K-3}
\end{gather}
We call $(c,u,w)\in \R\times\K\times\K$ satisfying (\ref{eq-tw-u}), (\ref{eq-tw-w})
a monotone (decreasing) travelling wave for \eqref{rd-u}, \eqref{rd-w}.\\
%The minimal speed of travelling waves is defined as
%\begin{eqnarray}\label{def-c0}
%       c_{\min} &:=& \inf\{c\in \R\,:\, \text{ there exists a monotone TW }(c,u,w)\}.
%\end{eqnarray}
\end{definition}
%The value of $c_{\min}=c_{\min}(k)$ in general depends on the parameter $k$ (and $\gamma$ as well).

We observe from \eqref{eq-tw-w} that the speed of a monotone (decreasing)
travelling
wave is always positive. 
\begin{theorem}[Existence of travelling waves]\label{the-exist}
For all $\gamma,k$ there is a positive number $c_{\min}<\infty$,
$c_{\min}=c_{\min}(\gamma,k)$, such 
that there 
exists a monotone travelling wave of \eqref{rd-u}, \eqref{rd-w} with speed $c$
for 
all $c\geq c_{\min}$ and such that there is no monotone travelling wave
with speed $c< c_{\min}$.
\end{theorem}
The value $c_{\min}$ gives thus the \emph{minimal speed} of travelling waves for
\eqref{rd-u}, \eqref{rd-w}.

In the remainder of this section we prove Theorem \ref{the-exist}.

For $0<\eps<(2\gamma)^{-1}$ we consider the
following strictly monotone and strictly parabolic approximation of
\eqref{rd-u}, \eqref{rd-w},
\begin{eqnarray}
  \partial_t \tilde{u}_\eps &=&
  \partial_{xx}\tilde{u}_\eps -\tilde{u}_\eps +\tilde{w}_\eps -\gamma k\tilde{u}_\eps(1-\tilde{w}_\eps),\label{rd-u-eps}\\
  \partial_t \tilde{w}_\eps &=&
  \eps\partial_{xx}\tilde{w}_\eps +\eps(\tilde{u}_\eps-\tilde{w}_\eps)+k\tilde{u}_\eps(1-\tilde{w}_\eps),\label{rd-w-eps}
\end{eqnarray}
and the corresponding travelling-wave equations
\begin{eqnarray}
  0 &=& u_\eps^{\prime\prime} +cu_\eps^\prime-u_\eps+w_\eps -\gamma k u_\eps(1-w_\eps)\label{tw-u-eps},\\
  0 &=& \eps w_\eps^{\prime\prime} +cw_\eps^\prime +\eps(u_\eps-w_\eps) +ku_\eps(1-w_\eps).\label{tw-w-eps}
\end{eqnarray}
Observe that the ODE system \eqref{tw-u-eps}, \eqref{tw-w-eps} has the same stationary points
as the original system \eqref{eq-tw-u}, \eqref{eq-tw-w}. For $\eps=0$ the systems \eqref{tw-u-eps}, \eqref{tw-w-eps}
and \eqref{eq-tw-u}, \eqref{eq-tw-w} coincide. The existence
of travelling waves for the auxiliary problem and a variational 
characterisation of the minimal speed follows from \cite{Vo}.
To state their result, we define functionals \mbox{$\Phi_1^\eps,\Phi_2^\eps: \K\times\K\to\R$},
\begin{eqnarray*}
  \Phi_1^\eps(\varrho,\sigma) &:=&\sup_{x\in\R}
  \frac{\varrho^{\prime\prime}(x)-\varrho(x) +\sigma(x)
        -\gamma k\varrho(x)(1-\sigma(x))}{-\varrho^\prime(x)},\\
  \Phi_2^\eps(\varrho,\sigma) &:=&
        \sup_{x\in\R}\,\frac{\eps\sigma^{\prime\prime}(x) +\eps(\varrho-\sigma)(x)+k\varrho(x)(1-\sigma(x))}
            {-\sigma^\prime(x)}.
\end{eqnarray*}
\begin{lemma}[\cite{Vo} Theorem I.4.2]\label{lem:volpert}
For all $c\geq c_\eps$, where $c_\eps\geq 0$ is defined by
\begin{eqnarray}
  c_\eps &:=& \inf_{\sigma,\varrho\in\K}
  \max\Big(\Phi_1^\eps(\varrho,\sigma),\Phi_2^\eps(\varrho,\sigma)\Big),
  \label{eq:def-c-eps} 
\end{eqnarray}
there exists a monotone travelling wave $(c,u_\eps,w_\eps)$ and $c_\eps$ is the minimal speed
of travelling waves for \eqref{rd-u-eps}, \eqref{rd-w-eps}. \hfill$\Box$
\end{lemma}
To prove the existence of travelling waves for \eqref{rd-u}, \eqref{rd-w} we pass in
\eqref{tw-u-eps}, \eqref{tw-w-eps} to the limit $\eps\to 0$.
First we derive bounds for $u_\eps,w_\eps$ which are uniform in $\eps>0$.
\begin{lemma}
Consider $\eps>0$ and let
$(c,u_\eps,w_\eps)$ be a monotone travelling wave for \eqref{rd-u-eps}, \eqref{rd-w-eps}.
Define $\mu,\nu>0$ to be
the positive solutions of
\begin{eqnarray}
  \mu^2+c\mu-1 &=& 0,\label{eq:def-mu}\\
  \nu^2-c\nu-(1+\gamma k) &=& 0.\label{eq:def-nu}
\end{eqnarray}
Then
\begin{eqnarray}
  \label{eq:bound-u-eps-prime}
  u_\eps^\prime &\geq& -\mu(1-u_\eps),\\
  u_\eps^\prime &\geq& -\nu u_\eps,
    \label{eq:bound-u-eps-prime+}\\
  |u_\eps^{\prime\prime}| &\leq& c\mu+1+\gamma k,\label{eq:bound-u-eps-prime2}\\
  \label{eq:bound-w-eps-prime}
  w_\eps^\prime &\geq& -\frac{\eps+k}{c}(1-w_\eps)
\end{eqnarray}
hold.
\end{lemma}
% In particular, for $c>\delta$ and any $l\in \N_0$ we obtain
%\begin{eqnarray}\label{eq:bounds-u,w-prime}
%  \|u_\eps\|_{C_l(\R)}, \|w_\eps\|_{C^l(\R)} &\leq& C(\gamma,k,\delta).
%\end{eqnarray}
\begin{proof}
We fix $x\in\R$ and define for $y\in (-\infty,x]$
\begin{eqnarray*}
  \hat{u}(y) &:=& 1-(1-u_\eps(x))e^{\mu (y-x)},\\
  \hat{w}(y) &:=& 1-(1-w_\eps(x))e^{L (y-x)},\quad\text{ where } L:=\frac{k+\eps}{c}.
\end{eqnarray*}
Then we observe that
\begin{gather*}
  \lim_{y\to -\infty}\hat{u}(y)\,=\,1,\quad \hat{u}(x)=u_\eps(x),\\
  \lim_{y\to -\infty}\hat{w}(y)\,=\,1,\quad \hat{w}(x)=w_\eps(x).
\end{gather*}
Moreover
\begin{gather*}
  \hat{u}^{\prime\prime} +c\hat{u}^\prime-\hat{u}+w_\eps-\gamma k\hat{u}(1-w_\eps)\\
  =\, -(1-u_\eps(x))\big(\mu^2+c\mu -1\big)e^{\mu (y-x)} -(1-w_\eps)-\gamma k\hat{u}(1-w_\eps)
  \,<\,0
\end{gather*}
and
\begin{gather*}
  \eps\hat{w}^{\prime\prime} +c\hat{w}^\prime-\eps\hat{w}+\eps u_\eps +ku_\eps(1-\hat{w})\\
  =\, -(1-w_\eps(x))\big(\eps L^2+cL-(\eps+k u_\eps)\big)e^{L (y-x)} -\eps(1-u_\eps)(x)\,<\, 0.
\end{gather*}
By the maximum principle, which we apply once for the scalar equation \eqref{tw-u-eps} and
once for the scalar equation \eqref{tw-w-eps}, we deduce that
\begin{gather*}
  u_\eps\,\leq\, \hat{u},\quad w_\eps\,\leq\, \hat{w}\qquad\text{ on }(-\infty,x],\label{eq:dom--u-eps}
\end{gather*}
thus
\begin{gather*}
  u_\eps^\prime(x)\,\geq \hat{u}^\prime(x)=-\mu\big(1-u_\eps(x)\big),\\
  w_\eps^\prime(x)\,\geq \hat{w}^\prime(x)=-L\big(1-w_\eps(x)\big),
\end{gather*}
which proves \eqref{eq:bound-u-eps-prime}, \eqref{eq:bound-w-eps-prime}.
The estimate \eqref{eq:bound-u-eps-prime2} follows from \eqref{tw-u-eps} and \eqref{eq:bound-u-eps-prime}.\\
In addition, by comparing $u_\eps$ on $[x,\infty)$ with
\begin{eqnarray*}
  \check{u}(y) &:=& u_\eps(x)e^{-\nu (y-x)}
\end{eqnarray*}
we obtain
\begin{eqnarray}
  \check{u}&\leq& u_\eps\quad\text{ on }[x,\infty),\label{eq:dom+-u-eps}\\
    \check{u}^\prime(x) &\leq& u_\eps^\prime(x),\notag
\end{eqnarray} 
which yields \eqref{eq:bound-u-eps-prime+}.
\end{proof}
%We now define%choose now a subsequence $\eps_i\to 0\,(i\to\infty)$ such that
%\begin{eqnarray}\label{eq:c-i}
%  c_0 &:=& \liminf_{\eps>0} c_\eps
%\end{eqnarray}
%and prove the existence of travelling waves for all speeds $c\geq c_0$.
By similar arguments one proves that corresponding properties
hold for travelling-wave solutions of \eqref{rd-u}, 
\eqref{rd-w}.
\begin{lemma}
Let $(c,u,w)$ be a monotone travelling wave for \eqref{rd-u}, \eqref{rd-w}.
Define $\mu,\nu>0$ to be
the positive solutions of
\begin{eqnarray}
  \mu^2+c\mu-1 &=& 0,\label{eq:def-mu-0}\\
  \nu^2-c\nu-(1+\gamma k) &=& 0.\label{eq:def-nu-0}
\end{eqnarray}
Then
\begin{eqnarray}
  \label{eq:bound-u-eps-prime-0}
  u^\prime &\geq& -\mu(1-u),\\
  u^\prime &\geq& -\nu u,
    \label{eq:bound-u-eps-prime+-0}\\
  |u^{\prime\prime}| &\leq& c\mu+1+\gamma
   k,\label{eq:bound-u-eps-prime2-0}\\ 
  \label{eq:bound-w-eps-prime-0}
  w^\prime &\geq& -\frac{k}{c}(1-w)
\end{eqnarray}
hold. \hfill$\Box$
\end{lemma}
We prove now the first statement in Theorem \ref{the-exist}.
\begin{proposition}\label{prop:ex-tw-c0}
For each $c\geq c_0$, where
\begin{eqnarray}\label{eq:c-i}
  c_0 &:= & \liminf_{\eps\to 0} c_\eps,
\end{eqnarray}
there exists a monotone travelling wave $(c,u,w)$ for \eqref{rd-u}, \eqref{rd-w}.
Moreover, the value $c_0$ is finite.
% Moreover,
%with $\mu,\nu$ as in \eqref{eq:def-mu}, \eqref{eq:def-nu}
%\begin{eqnarray}
% \label{eq:bound-u-prime}
%  u^\prime &\geq& -\mu(1-u),\\
%  u^\prime &\geq& -\nu u, \label{eq:bound-u-prime+}\\
%  |u^{\prime\prime}| &\leq& c\mu+1+\gamma k,\label{eq:bound-u-prime2}\\
%  \label{eq:bound-w-prime}
%  w^\prime &\geq& -\frac{k}{c}(1-w)
%\end{eqnarray}
%holds.
\end{proposition}
\begin{proof}
Assume first that $c_0<\infty$ and fix an arbitrary $c\geq c_0$ and a
subsequence $\eps_i\to 0\, (i\to\infty)$ such that 
\begin{eqnarray*}
  c_0 &=& \lim_{i\to\infty} c_{\eps_i}.
\end{eqnarray*}
By Lemma \ref{lem:volpert} there exists a sequence of monotone travelling waves
$(c_i,u_{\eps_i},w_{\eps_i})$ for \eqref{rd-u-eps}, 
\eqref{rd-w-eps}, with $c_i=c$ if $c>c_0$ and $c_i=c_{\eps_i}$ if
$c=c_0$ such that 
\begin{eqnarray}
  0 &=& u_{\eps_i}^{\prime\prime} +c_iu_{\eps_i}^\prime-u_{\eps_i}+w_{\eps_i}
  -\gamma k u_{\eps_i}(1-w_{\eps_i}),\label{eq:tw-u-i}\\
  0 &=& \eps_i w_{\eps_i}^{\prime\prime} +c_iw_{\eps_i}^\prime +\eps_i(u_{\eps_i}-w_{\eps_i})
  +ku_{\eps_i}(1-w_{\eps_i}).\label{eq:tw-w-i}
\end{eqnarray}
Since travelling waves are invariant under space shifts, we can assume
without loss of generality that
\begin{gather}
  u_{\eps_i}(0)=\frac{1}{2}\quad\text{ for all }i\in\N.\label{eq:u0=0.5}
\end{gather}
By $0< u_{\eps_i},w_{\eps_i}<1$, the monotonicity of $w_{\eps_i}$ and
\eqref{eq:bound-u-eps-prime}, \eqref{eq:bound-u-eps-prime2} 
there exists $u\in C^{1,1}(\R), w\in L^\infty(\R)$ such that for all
$0<\alpha<1$, $R>0$
\begin{eqnarray}
  u_{\eps_i} &\to& u\quad\text{ in }C^{1,\alpha}([-R,R]),\label{eq:conv-u-i}\\
  w_{\eps_i} &\to& w\quad\text{ pointwise almost everywhere in }\R\label{eq:conv-w-i}
\end{eqnarray}
hold for a subsequence $\eps_i\to 0 (i\to\infty)$.

Multiplying \eqref{eq:tw-u-i}, \eqref{eq:tw-w-i} by a function $\eta\in C^\infty_c(\R)$ and
integrating we deduce
\begin{eqnarray*}
  0 &=& \int_\R \Big(- \eta^\prime u_{\eps_i}^{\prime} +c_{i}\eta
  u_{\eps_i}^\prime-\eta(u_{\eps_i}-w_{\eps_i}) -\eta 
  \gamma k u_{\eps_i}(1-w_{\eps_i})\Big),\\
  0 &=& \int_\R\Big( \eps_i \eta^{\prime\prime}w_{\eps_i}
  -c_{i}\eta^\prime w_{\eps_i} +{\eps_i}\eta (u_{\eps_i}-w_{\eps_i}) 
  +k\eta u_{\eps_i}(1-w_{\eps_i})\Big).
\end{eqnarray*}
Due to \eqref{eq:conv-u-i}, \eqref{eq:conv-w-i} we can pass  to the
limit $\eps_i\to 0$ in these equations
and get
\begin{eqnarray*}
  0 &=& \int_\R \Big(- \eta^\prime u^{\prime} +c\eta u^\prime-\eta(u-w) -\eta
  \gamma k u(1-w)\Big),\\
  0 &=& \int_\R\big(  -c\eta^\prime w +k\eta u(1-w)\big).
\end{eqnarray*}
It follows that $(c,u,w)$ solve \eqref{eq-tw-u}, \eqref{eq-tw-w} and, by
a bootstrapping argument, that $u,w$ are smooth.
Moreover, \eqref{eq:u0=0.5} and \eqref{eq:conv-u-i} yield that $u(0)=1/2$
and by Lemma 
\ref{lem-tw-basic} we obtain that $c>0$ and that $(c,u,w)$ 
is a monotone travelling wave.
% By \eqref{eq:bound-u-eps-prime}, \eqref{eq:bound-u-eps-prime+}
%we deduce \eqref{eq:bound-u-prime}, \eqref{eq:bound-u-prime+}
%and from \eqref{eq:bound-w-eps-prime} we get \eqref{eq:bound-w-prime}.
Since $c\geq c_0$ was arbitrary, the first part of the Proposition is
proved. To prove that $c_0<\infty$, let $\varrho,\sigma$ be two smooth
strictly monotonically decreasing functions with
\begin{alignat*}{2}
  \varrho(x)\,&=\,\sigma(x)\,=\,e^{-x}\quad&&\text{ for all }x\geq 1,\\
  \varrho(x)\,&=\,\sigma(x)\,=\,1-e^{-x}\quad&&\text{ for all }x\leq -1.
\end{alignat*}
Then there exists a constant $C(k)<\infty$ such that
for all $0<\eps<1$ and all $|x|\geq 1$ 
\begin{align*}
  \frac{\varrho^{\prime\prime}(x)-\varrho(x) +\sigma(x)
    -\gamma k\varrho(x)(1-\sigma(x))}{-\varrho^\prime(x)}\,&\leq\, C(k),\\
  \frac{\eps\sigma^{\prime\prime}(x) +\eps(\varrho-\sigma)(x)+k\varrho(x)(1-\sigma(x))}
            {-\sigma^\prime(x)}\,&\leq\, C(k).
\end{align*}
We can estimate the same ratios for all $0<\eps<1$ and $|x|\leq 1$ by
a constant depending only on $\|\varrho\|_{C^2([-1,1])}$,
$\|\sigma\|_{C^2([-1,1])}$. By Lemma \ref{lem:volpert} and the
definition of $c_\eps$ in \eqref{eq:def-c-eps} it follows that
\begin{gather*}
  c_\eps\,\leq\, C(k)\quad\text{ for all }0<\eps<1.
\end{gather*}
In particular, $c_0\leq C(k)$ is finite.
\end{proof}
We now complete the proof of Theorem \ref{the-exist}.
\begin{proposition}
\label{prop:c0=cmin}
There is no monotone travelling wave for \eqref{rd-u}, \eqref{rd-w} with speed $c<c_0$.
\end{proposition}
\begin{proof}
Assume that $(c,u,w)$ with $c<c_0$, $u,w\in\K$ satisfies \eqref{eq-tw-u}, \eqref{eq-tw-w}.
We then obtain that
\begin{align}
  \Phi^\eps_1(u,w)\,&=\,c,\notag\\
  \Phi^\eps_2(u,w)\,&\leq\, c +
  \eps\sup_{x\in\R}\frac{w^{\prime\prime}(x) + u(x)-w(x)}{-w^\prime(x)}.\label{eq:est-phi2}
\end{align} 
From \eqref{eq-tw-w} and $u<1$ we deduce that
\begin{gather}
        \frac{u-w}{-w^\prime}\,\leq\, \frac{u(1-w)}{-w^\prime}\,=\,\frac{c}{k}.\label{eq:phi2-term1}
\end{gather}
A differentiation in \eqref{eq-tw-w} yields that
\begin{align*}
        -cw^{\prime\prime}\,&=\, ku^\prime(1-w) -kuw^\prime
                \,=\, -cw^\prime\frac{u^\prime}{u} -kuw^\prime
\end{align*}
which gives, together with \eqref{eq:bound-u-eps-prime+-0}, that
\begin{gather}
 \frac{w^{\prime\prime}}{-w^\prime}\,=\, \frac{-u^\prime}{u} -\frac{k}{c}u\,\leq\ \nu
 \label{eq:phi2-term2}
\end{gather}
with $\nu$ as in \eqref{eq:def-nu-0}.
Using \eqref{eq:phi2-term1}, \eqref{eq:phi2-term2} in \eqref{eq:est-phi2} we obtain that
\begin{gather*}
 \Phi^\eps_2(u,w)\,\leq\ c+\eps C,
\end{gather*}
where $C$ is independent of $\eps>0$.
From Lemma \ref{lem:volpert} and the definition of $c_0$,
\eqref{eq:def-c-eps} and \eqref{eq:c-i} we deduce that
\begin{gather*}
 c_0\,\leq\, \liminf_{\eps\to 0} \max (\Phi^\eps_1(u,w),\Phi^\eps_2(u,w))\,=\, c,
\end{gather*}
which is a contradiction to our assumption $c<c_0$.
\end{proof}

%===================================================
% section limit
%===================================================
\section{The Stefan-like limit problem}
\label{sec-limit}
The reaction-diffusion system \eqref{rd-u}, \eqref{rd-w} converges to a
Stefan-like free boundary problem as $k$ tends to infinity, see \cite{HKR}.  
For solutions ${U}_\infty,{W}_\infty$ of this limit problem
\begin{eqnarray}
  0 &=& U_\infty(1-W_\infty),\label{st-w0}
\end{eqnarray}
holds and the spatial
domain splits in a region where ${U}_\infty=0$, and a region where
${U}_\infty>0$ and ${W}_\infty=1$. If we denote 
their common boundary at time $t$ by $\Gamma(t)$ then $U_\infty,
W_\infty$ satisfy
\begin{alignat}{2}
  \partial_t {U}_\infty \,&=\, \Delta {U}_\infty
  -{U}_\infty + 1&\quad&\text{ in }\{{W}_\infty=1\},\label{st-w=1}\\
  \gamma \partial_t {W}_\infty \,&=\, {W}_\infty&\quad&\text{ in
  }\{{U}_\infty=0\}\label{st-u=0}
\end{alignat}
and a continuity and jump condition on $\Gamma(t)$,
\begin{align}
  \mbox{}[U_\infty(t,.)]\,&=\, 0,\label{st-Gamma2}\\
   - [\nabla{U}_\infty(t,.)\cdot\nu(t,.)]\,&=\,\gamma
  [{W}_\infty(t,.)]\vec{v}(t,.)\cdot\nu(t,.)\label{st-Gamma}, 
\end{align}
where $\vec{v}(t,.)$ and $\nu(t,.)$ are the velocity and the unit normal of the free boundary
$\Gamma(t)$, pointing into
$\{U_\infty=0\}$, and $[.]$ denotes the jump across the free boundary from
the region $\{U_\infty>0\}$ to $\{U_\infty=0\}$.\\
As for the reaction-diffusion system, travelling-wave solutions
$U_\infty, W_\infty$ are
given by a speed $c\in\R$ and \emph{profile functions} 
$u_\infty,w_\infty$,
\begin{gather*}
  U_\infty(x,t)\,=\, u_\infty(x-ct),\quad W_\infty(x,t)\,=\,
 w_\infty(x-ct).
\end{gather*}
We are interested in monotone travelling waves which
connect unity and zero,
\begin{align}
  u_\infty(x),w_\infty(x) &\,\to\, 1\quad\text{ as }x\to -\infty,\label{eq:tw-lim-infty-}\\
  u_\infty(x),w_\infty(x) &\,\to\, 0\quad\text{ as }x\to \infty,\label{eq:tw-lim-infty}\\
  u_\infty,w_\infty \quad&\text{ are monotonically decreasing}.\label{eq:tw-lim-mono}
\end{align}
Due to the shift invariance of travelling waves, the condition
\eqref{st-w0} and the monotonicity of $u_\infty, w_\infty$ we can assume
that 
\begin{alignat}{2}
  u_\infty(x)\,&=\,0&\qquad&\text{ for }x\geq 0,\label{eq:tw-st0}\\
  u_\infty(x)\,&>\,0,\, w_\infty(x)\,=\,1&\quad&\text{ for }x<0.\label{eq:tw-st1}
\end{alignat}
From \eqref{st-w=1}-\eqref{st-Gamma} we then obtain
\begin{alignat}{2}
  0\,&=\, u_\infty^{\prime\prime}+cu_\infty^\prime -u_\infty +1&\quad&\text{
    for }x<0,\label{eq:tw-st2}\\
  0\,&=\, \gamma c w_\infty^\prime + w_\infty&\quad&\text{ for }x>0,
    \label{eq:tw-st3}
\end{alignat}
and the continuity and jump condition
\begin{align}
      u_\infty(0-)\,&=\,0,\label{eq:tw-st5}\\
  \gamma c (1-w_\infty(0+)) &\,=\,
    -u^\prime_\infty(0-)\label{eq:tw-st4}.
\end{align}
\begin{proposition}
\label{prop:ex-tw-lim}
For all $c\geq c_\infty$, where
\begin{eqnarray}\label{eq-c-infty}
        c_\infty &:=& \frac{1}{\sqrt{\gamma (1+\gamma)}},
\end{eqnarray}
there exists a unique solution $u_\infty,w_\infty$ of
\eqref{eq:tw-lim-infty-}-\eqref{eq:tw-st4}.
This solution is given by
\begin{eqnarray}
        u_\infty(x) &=&
        \begin{cases}
                1- e^{\alpha x} &\text{ if } x\leq 0,\\
                0 &\text{ if }x\geq 0,
        \end{cases}
        \label{tw-lim-u}\\
        w_\infty(x) &=&
        \begin{cases}
                1 &\text{ if }x<0,\\
                \beta e^{\textstyle{-\frac{x}{c\gamma}}} &\text{ if }x>0,
        \end{cases}\label{tw-lim-w}
\end{eqnarray}
where
\begin{gather}
        \alpha=\alpha(c) := -\frac{c}{2} + \sqrt{\frac{c^2}{4}+1},\label{eq:def-alpha-mu}\\
        \qquad \beta=\beta(c):=1-\frac{\alpha(c)}{\gamma c}.\label{eq:def-alpha}
\end{gather}
For $c < c_\infty$, there does not exist any solution.
\end{proposition}
\begin{proof}
We deduce that \eqref{eq:tw-st2} and \eqref{eq:tw-lim-infty-},
\eqref{eq:tw-st5} hold if and only if
\eqref{tw-lim-u} holds, with $\alpha$ as in \eqref{eq:def-alpha-mu}.
Similarly \eqref{eq:tw-st3} and \eqref{eq:tw-lim-infty} are satisfied if
and only if $w_\infty$ satisfies \eqref{tw-lim-w} with $\beta\in\R$.
Since
\begin{gather*}
  1- w_\infty(0+)\,=\, 1-\beta,\quad
  -u_\infty^\prime(0-)\,=\, \alpha
\end{gather*}
we obtain that the jump condition \eqref{eq:tw-st4} is satisfied if and
only if $\beta$ satisfies \eqref{eq:def-alpha}.
Finally, by \eqref{tw-lim-w} and \eqref{eq:def-alpha}, the solution
$w_\infty$ is nonnegative if 
and only if $c\geq c_\infty$.
\end{proof}
%===================================================
% section speed selection
%===================================================
\section{Selection of the minimal speed}
\label{sec-speed}
In this section we prove that for sufficiently large values $k>k_0$ a nonlinear selection principle determines
the minimal speed of travelling waves. The threshold $k_0$ is obtained explicitly in terms of
the given constant $\gamma$.
First, we have to analyse the behaviour of
travelling-wave solutions at infinity.
\subsection{Behaviour at infinity}
The linear selection principle for the minimal speed is based on the analysis
of the linearised system at the unstable stationary point $(u,w)=(0,0)$ of \eqref{eq-tw-u}, \eqref{eq-tw-w}.
In the next lemma we show that solutions have to decay exponentially to zero as $x$ tends to infinity.
\begin{lemma}\label{lem-decay-infty+}
Let $(c,u,w)$ be a monotone travelling wave. Then
\begin{gather}\label{eq-decay-infty+}
        \lim_{x\to\infty}\frac{1}{x}\log u(x)\,=\,\lim_{x\to\infty}\frac{1}{x}\log w(x)\,=\,\lambda,
\end{gather}
where $\lambda$ is a negative root of the cubic equation
\begin{eqnarray}\label{eq-lambda}
        0 &=& \lambda^3 + c\lambda^2 -(1+\gamma k)\lambda -\frac{k}{c}.
\end{eqnarray}
\end{lemma}
\begin{proof}
With the definitions
\begin{gather*}
  V_0(x)\,:=\,
  \begin{pmatrix}
    u^\prime(x)\\
    u(x)\\
    w(x)
  \end{pmatrix}\!,
  \quad
  A_0\,:=\,
  \begin{pmatrix}
    -c & 1+\gamma k & -1\\
    1 &0 & 0\\
    0 &-\frac{k}{c} &0
  \end{pmatrix}\!,
  \quad
  B_0\,:=\,
  \begin{pmatrix}
    -\gamma\\
    0\\
    \frac{1}{c}
  \end{pmatrix}
\end{gather*}
the system \eqref{eq-tw-u}, \eqref{eq-tw-w} is equivalent to
\begin{eqnarray}
  V_0^\prime(x) &=& A_0V_0(x) + ku(x)w(x)B_0
\end{eqnarray}
and the linearized system at $(0,0,0)$ is given by
\begin{eqnarray}\label{eq:xi-lin}
  \xi^\prime(x) &=& A_0\xi(x).
\end{eqnarray}
The eigenvalues of $A_0$ are the
solutions of \eqref{eq-lambda}. Since $s\mapsto s^3 + c s^2 -(1+\gamma k)s -\frac{k}{c}$ is negative
at $s=0$ and becomes positive as $s\to +\infty$ there exists a positive eigenvalue $\lambda_+>0$.
The other two solutions of \eqref{eq-lambda} satisfy the equation
\begin{eqnarray}\label{eq-lambda-+}
  0 &=& \lambda^2 +(c+\lambda_+)\lambda +\frac{k}{c\lambda_+} 
\end{eqnarray}
which has, depending on the values of $c,k,\lambda_+$, either two negative roots, one repeated negative root
or two complex-conjugate roots with negative real part.
Since $V_0(x)$ converges to zero as $x\to +\infty$, the curve
$x\mapsto V_0(x)$ is for sufficiently large values of $x$ contained in the stable manifold and
converges exponentially to zero; see for example \cite{Pe} Section 2.7.
By \cite{CL} Theorem XIII.4.5 there exists a solution $\xi=(\xi_1,\xi_2,\xi_3)$ of the linearized system
\eqref{eq:xi-lin} and $\delta>0$ with
\begin{eqnarray}\label{eq:exp-close}
  V_0(x) &=& \xi(x) + O(e^{(\alpha-\delta)x}),\\
  \xi(x) &=& O(e^{\alpha x})\notag
\end{eqnarray}
as $x\to +\infty$, where $\alpha<0$ is the real part of an eigenvalue $\lambda$ of $A_0$.
One checks that $A_0$ has no eigenvector with a component equal to zero; therefore
\begin{eqnarray}\label{eq:dec-comp}
  \xi_i(x) &=& O(e^{\alpha x})\quad\text{ for } i=1,2,3
\end{eqnarray}
holds as $x\to +\infty$. Let us show that in fact $\lambda$ is real. Assume that
$\lambda=\alpha+i\beta$ with $\beta\neq 0$. Then
$\xi(x)$ describes, as $x\to +\infty$, a spiral around the origin contained
in the plane $P$ spanned by the real and imaginary part of an eigenvector
of $A_0$ with eigenvalue $\lambda$.
But then, since the difference between $\xi$ and $V_0$ decays exponentially faster than $\xi$,
$V_0$ has to take values outside the set $\{z=(z_1,z_2,z_3):z_1<0,z_2>0,z_3>0\}$,
which is a contradiction to the assumption that $(c,u,w)$ is a
monotone travelling wave. This shows that $\lambda=\alpha<0$ is real. Thus \eqref{eq-decay-infty+}
follows from \eqref{eq:exp-close} and \eqref{eq:dec-comp}.
\end{proof}
The equation \eqref{eq-lambda} connects the speed and the decay
rate at $+\infty$ of a travelling wave. We now further analyze
this relation.
\begin{lemma}\label{lem:curve-lambda-c}
For all $k>0$ and all $\lambda<0$ there exists a unique value
$\bar{c}_k(\lambda)>0$ such that $\lambda$ satisfies 
\eqref{eq-lambda} for $c=\bar{c}_k(\lambda)$. The function
$\lambda\mapsto \bar{c}_k(\lambda)$ 
attains a positive minimum $c_{\lin}=c_{\lin}(k)$ at a unique value
$\lambda_{\lin}=\lambda_{\lin}(k)$ and 
$c_{\lin},\lambda_{\lin}$ are given by
\begin{eqnarray}
  \lambda_{\lin}&<& 0\notag\\
  \label{eq-lambda-lin1}
  3\lambda_{\lin}^2 &=& -(1+\gamma k)+2\sqrt{(1+\gamma k)^2+3k},\\
  \lambda_{\lin} c_{\lin} &=& (1+\gamma k) -\sqrt{(1+\gamma k)^2+3k}.\label{eq-lambda-lin}
\end{eqnarray}
Moreover
\begin{alignat}{2}
  \bar{c}_k^\prime(\lambda)\,&<\, 0 \quad &&\text{ for
 }\lambda<\lambda_{\lin},\label{eq:mon-ck-1}\\
  \bar{c}_k^\prime(\lambda)\,&>\, 0 \quad &&\text{ for
 }\lambda>\lambda_{\lin},\label{eq:mon-ck-2}
\end{alignat}
with
\begin{eqnarray*}
  \bar{c}_k(\lambda) &\to& \infty\quad\text{ as }\lambda\to -\infty,\\
  \bar{c}_k(\lambda) &\to& \infty\quad\text{ as }\lambda\to 0.
\end{eqnarray*}
\end{lemma}
\begin{proof}
From \eqref{eq-lambda} we obtain that $\bar{c}_k$ is given by
\begin{eqnarray*}
  \bar{c}_k(\lambda) &=& -\frac{1}{2}\Big(\lambda-(1+\gamma k)\frac{1}{\lambda}\Big)
  +\sqrt{\frac{1}{4}\Big(\lambda-(1+\gamma k)\frac{1}{\lambda}\Big)^2+\frac{k}{\lambda^2}},
\end{eqnarray*}
that $\bar{c}_k$ is strictly positive and that $\bar{c}_k$ tends to
infinity as $\lambda\to -\infty$ or $\lambda\to 0$. 
Therefore the positive minimum $c_{\lin}$ is attained at a value
$\lambda_{\lin}$ and 
$\bar{c}^\prime(\lambda_{\lin})=0$ holds. By \eqref{eq-lambda} this implies
\begin{eqnarray}\label{eq:lam1}
  0 &=& 3\lambda_{\lin}^2 + 2\lambda_{\lin} c_{\lin} - (1+\gamma k).
\end{eqnarray}
Moreover, by \eqref{eq-lambda},
\begin{eqnarray}\label{eq:lam2}
        0 &=& \lambda_{\lin}^3 + c_{\lin}\lambda_{\lin}^2 -(1+\gamma
        k)\lambda_{\lin} -\frac{k}{c_{\lin}}. 
\end{eqnarray}
One checks that \eqref{eq-lambda-lin1}, \eqref{eq-lambda-lin} is equivalent
to \eqref{eq:lam1}, \eqref{eq:lam2}.
In particular, $\bar{c}_k^\prime$ has only one zero and we deduce that
$\bar{c}_k^\prime(\lambda)<0$ for $\lambda<\lambda_{\lin}$
and  $\bar{c}_k^\prime(\lambda)>0$ for $0>\lambda>\lambda_{\lin}$.
\end{proof}
\begin{corollary}\label{cor-c(lambda)}
The minimal speed $c_{\min}(k)$ of travelling waves satisfies the estimate
\begin{eqnarray}\label{eq-lower-bound}
  c_{\min}(k) &\geq& c_{\lin}(k).
\end{eqnarray}
\end{corollary}
\begin{proof}
According to Lemma \ref{lem-decay-infty+} for a monotone travelling wave with speed $c$, a negative root
$\lambda$ of \eqref{eq-lambda} exists.
On the other hand $c_{\lin}(k)$ gives the minimal value of $c$ such that \eqref{eq-lambda} has a negative solution.
\end{proof}
We now state the result corresponding to Lemma \ref{lem-decay-infty+}
if we consider $x$ approaching $-\infty$.
\begin{lemma}\label{lem-decay-infty-}
Let $(c,u,w)$ be a monotone travelling wave and let $k\geq 1$.
Then
\begin{eqnarray}
        \lim_{x\to-\infty}\frac{1}{x}\log \big(1-u(x)\big) &=& \mu,\label{eq-decay-infty-u}\\
        \lim_{x\to-\infty}\frac{1}{x}\log \big(1-w(x)\big) &=&\frac{k}{c},\label{eq-decay-infty-w}
\end{eqnarray}
where $\mu>0$ is given by
\begin{eqnarray}\label{eq-lambda-}
        \mu &=& -\frac{c}{2} +\sqrt{\frac{c^2}{4}+1}.
\end{eqnarray}
For $0<k<1$ the limits in \eqref{eq-decay-infty-u}, \eqref{eq-decay-infty-w}
are either equal to $\mu$ or equal to $k/c$.
\end{lemma}
\begin{proof}
With
\begin{gather*}
  V_1(x)\,:=\,
  \begin{pmatrix}
    -u^\prime(x)\\
    1-u(x)\\
    1-w(x)
  \end{pmatrix}\!,
  \quad
  A_1\,:=\,
  \begin{pmatrix}
    -c & 1 & -(1+\gamma k)\\
    1 &0 & 0\\
    0 & 0 &\frac{k}{c}
  \end{pmatrix}\!,
  \quad
  B_1\,:=\,
  \begin{pmatrix}
    \gamma\\
    0\\
    -\frac{1}{c}
  \end{pmatrix}
\end{gather*}
the system \eqref{eq-tw-u}, \eqref{eq-tw-w} is equivalent to
\begin{eqnarray}\label{eq:V-A-}
  V_1^\prime(x) &=& A_1V_1(x) + k\big(1-u(x)\big)\big(1-w(x)\big)B_1
\end{eqnarray}
and the linearized system at $V_1=(0,0,0)$ is given by
\begin{eqnarray}\label{eq:xi-lin-}
  \zeta^\prime(x) &=& A_1\zeta(x).
\end{eqnarray}
The matrix $A_1$ has the positive eigenvalues
\begin{gather*}
  \mu\,=\, -\frac{c}{2} +\sqrt{\frac{c^2}{4}+1},\quad
 \mu_1\,=\,\frac{k}{c}
\end{gather*}
and
one negative eigenvalue 
\begin{gather*}
  \mu_2\, =\, -\frac{c}{2} - \sqrt{\frac{c^2}{4} +1}.
\end{gather*}
Since $V_1(x)\to 0$ as $x\to -\infty$ we deduce
that $V_1(x)$ is for sufficiently small $x$ contained in the unstable manifold
of \eqref{eq:V-A-} at $(0,0,0)$.
 Using \cite{CL} Theorem XIII.4.5 we obtain the existence of a solution
$\zeta$ of the linearized system
\eqref{eq:xi-lin-} and a $\delta>0$ with
\begin{eqnarray}
  V_1(x) &=& \zeta(x) + O(e^{(\beta+\delta)x}),\label{eq:lin-close-}\\
  \zeta(x) &=& O(e^{\beta x})\label{eq:lin-decay-}
\end{eqnarray}
as $x\to -\infty$, where $\beta=k/c$ or $\beta=\mu$.
One checks that $k/c>\mu$ holds for $k\geq 1$.
Thus, if $\beta=k/c$ then the trajectory of $V_1$ as $x\to -\infty$ has to be tangential
to the eigenspace of $A_1$ corresponding to the eigenvalue $k/c$. On the other hand this
eigenspace is spanned by a vector $(k/c,1,a)$, where $a<0$ for $k\geq 1$.
Therefore $V_1$ leaves the region $\{z=(z_1,z_2,z_3): z_2,z_3>0\}$, which
is a contradiction to $u(x),w(x)<1$.
This proves that $\beta=\mu$ for $k\geq 1$ and implies that the trajectory of
$\zeta(x)$ for $x\to -\infty$ is tangential to the eigenspace corresponding
to the eigenvalue $\mu$, which is spanned by the vector $(\mu,1,0)$.
Therefore
\begin{eqnarray*}
  \zeta_2(x) &=& O(e^{\mu x})
\end{eqnarray*}
and \eqref{eq-decay-infty-u} follows from \eqref{eq:lin-close-}. By
\eqref{eq-tw-w} we deduce
\begin{eqnarray*}
  \frac{d}{dx}\log \big(1-w(x)\big) &=& \frac{k}{c} + O(e^{\mu x})\quad\text{ as }x\to -\infty,
\end{eqnarray*}
and thus \eqref{eq-decay-infty-w} holds.
\end{proof}
\begin{corollary}\label{cor:decay-}
Let $k\geq 1$. If we consider two different monotone travelling waves
for \eqref{rd-u}, \eqref{rd-w}, 
the one with the lower speed converges faster to unity as $x$ approaches $-\infty$.
\end{corollary}
\begin{proof}
For a monotone travelling wave $(c,u,w)$, by Lemma \ref{lem-decay-infty-} the convergence of $u,w$
to unity as $x\to -\infty$
is exponential with rate $\mu$ for $u$ and rate $k/c$ for $w$, where $\mu$ is given by \eqref{eq-lambda-}. Since
\begin{gather*}
  \frac{d}{dc}\Big( -\frac{c}{2} +\sqrt{\frac{c^2}{4}+1}\Big)\,=\,
  \frac{1}{2}\Big(-1+\frac{c}{\sqrt{c^2+4}}\Big)\,<\, 0
\end{gather*}
we deduce that both convergence rates are decreasing with $c$. 
\end{proof}
\subsection{The reduced system for $c=c_\infty$}
It was observed in \cite{Ki} that the travelling-wave equations \eqref{eq-tw-u}, \eqref{eq-tw-w}
are for the speed $c=c_\infty$ defined in
\eqref{eq-c-infty} remarkable in being equivalent to a system of two first-order equations.
\begin{lemma}\label{lem:red}
Let $u\in \K_0,w\in \K$. Then $(c_\infty,u,w)$ is a monotone travelling wave if and only if $u,w$ satisfy
\begin{eqnarray}\label{eq-red-u}
        u^\prime &=& \gamma c_\infty (u-w),\\
        w^\prime &=& -k\gamma (1+\gamma) c_\infty u(1-w).\label{eq-red-w}
\end{eqnarray}
\end{lemma}
\begin{proof}
Multiplying \eqref{eq-red-w} by $c_\infty$ and using \eqref{eq-c-infty}
we see that \eqref{eq-red-w} and \eqref{eq-tw-w} are equivalent for $c=c_\infty$.
Next we obtain from \eqref{eq-c-infty} that
\begin{eqnarray}
  \big(u^\prime -\gamma c_\infty(u-w)\big)^\prime &=& u^{\prime\prime} +c_\infty u^\prime
  -\frac{1}{\gamma c_\infty}u^\prime -\gamma k u(1-w)\label{eq:aux-red}
\end{eqnarray}
and we see that \eqref{eq-red-u} implies \eqref{eq-tw-u}.
Conversely, from \eqref{eq-tw-u} and \eqref{eq:aux-red} we deduce that
\begin{eqnarray*}
  \big(u^\prime -\gamma c_\infty(u-w)\big)^\prime &=& -\frac{1}{\gamma c_\infty}\big(u^\prime-\gamma c_\infty(u-w)\big),
\end{eqnarray*}
whose solutions are
\begin{eqnarray*}
  \big(u^\prime -\gamma c_\infty(u-w)\big)(x) &=& a e^{-\frac{1}{\gamma
  c_\infty}x}\quad\text{ for any } a\in\R. 
\end{eqnarray*}
The condition that $u(x),w(x)$ converge exponentially to unity as $x\to -\infty$
implies that $a=0$ and \eqref{eq-red-u} therefore holds.
\end{proof}
This reduction allows us to prove the existence of a travelling wave with speed $c_\infty$
for \eqref{rd-u}, \eqref{rd-w} by a phase-plane analysis for
\eqref{eq-red-u}, \eqref{eq-red-w}.
\begin{proposition}\label{prop-red}
For all $\gamma,k>0$ there exists a monotone travelling wave  $(c_\infty,u_0,w_0)$
for \eqref{rd-u}, \eqref{rd-w}. As $x$ tends to infinity, $u_0,w_0$ decay exponentially to
zero with decay rate
\begin{gather}\label{eq:dec-u-0}
  \lambda_{\infty}\,=\,\lambda_\infty(k) \,:=\,
  c_\infty\gamma\Big(\frac{1}{2}-\sqrt{\frac{1}{4}+k(1+\gamma)}\Big)\,<\,0. 
\end{gather}
\end{proposition}
\begin{proof}
The system \eqref{eq-red-u}, \eqref{eq-red-w}
has the two stationary points $(0,0)$ and $(1,1)$.
We define the set
\begin{gather}\label{def:J}
  J\,:=\,\{(x_1,x_2)^T\in\R^2 : 0\leq x_1<x_2\leq 1\}
\end{gather}
and observe that $\R^2\setminus J$ is an invariant region for
\eqref{eq-red-u}, \eqref{eq-red-w} (see Figure \ref{fig:red}).\\
Next we consider the linearisation at $(0,0)$ of \eqref{eq-red-u}, \eqref{eq-red-w}
which is given by
\begin{eqnarray*}
  \zeta^\prime &=& c_\infty\gamma
  \begin{pmatrix}
    1 & -1\\
    -k(1+\gamma) & 0
  \end{pmatrix}
  \zeta,
\end{eqnarray*}
for $\xi:\R\to\R^2$.
The eigenvalues of this linear system are $\lambda_\infty$, as defined
in \eqref{eq:dec-u-0}, and 
\begin{gather*}
  c_\infty\gamma\Big(\frac{1}{2}+\sqrt{\frac{1}{4}+k(1+\gamma)}\Big)>0.
\end{gather*}
An eigenvector with eigenvalue $\lambda_\infty$ is given by
\begin{gather*}
  v \,:=\,
  \begin{pmatrix}
    -\frac{1}{2} + \sqrt{\frac{1}{4}+k(1+\gamma)}\\
    k(1+\gamma)
  \end{pmatrix}
\end{gather*}
One checks that $v_1<v_2$ holds for the components of $v$ and deduces
that the eigenspace corresponding to $\lambda_\infty$ intersects the
set $J$ defined in \eqref{def:J}.
\begin{figure}[t]
\centering
\includegraphics[height=65mm]{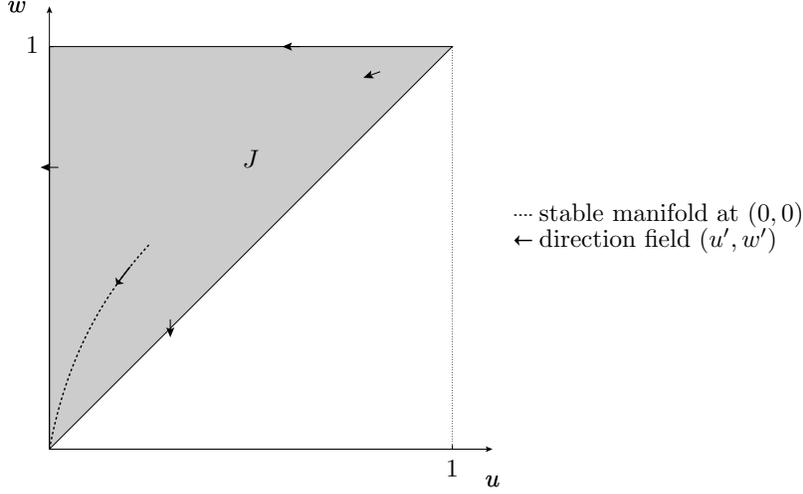}
\caption{Schematic of the phase portrait of
  \eqref{eq-red-u},\eqref{eq-red-w}. The dashed area represents the set
  $J$, the incoming arc the stable manifold at $(0,0)$ and the arrows the
  direction field.}
\label{fig:red}
\end{figure}

By the stable manifold theorem, see for example \cite{Pe}, Theorem 2.7, there exists a
trajectory that, taking $t$ as its parameter variable, converges to $(0,0)$ as $t\to\infty$
and starts at $t=0$ in $J$, since the stable manifold is tangential to the eigenspace corresponding
to the eigenvalue $\lambda_\infty$.
Following this trajectory back with decreasing $t$ we cannot leave $J$, since otherwise
the trajectory would stay in $\R^2\setminus J$ as $t$ increases and thus
could not reach any point in $J$ at $t=0$. In $J$, with decreasing $t$, the trajectory
has to be monotone in both components and therefore has to approach the
stationary point $(1,1)$. Thus the trajectory connects $(1,1)$ to $(0,0)$ and satisfies
\eqref{eq-red-u}, \eqref{eq-red-w}. By Lemma \ref{lem:red}
this shows that $(c_\infty,u_0,w_0)$ is a monotone travelling wave for
\eqref{rd-u}, \eqref{rd-w}. Moreover $u_0(x),w_0(x)$ are, for sufficiently large $x$, in the stable
manifold of \eqref{eq-red-u}, \eqref{eq-red-w}
and we deduce that they converge exponentially fast to zero, with
decay rate given by $\lambda_\infty$.   
\end{proof}
The existence of a travelling wave with speed $c_\infty$ implies immediately
the following estimate.
\begin{corollary}
The minimal speed $c_{\min}(k)$ of travelling waves satisfies
\begin{eqnarray}\label{eq-upper-bound}
  c_{\min}(k) &\leq& c_\infty.
\end{eqnarray}
\end{corollary}
As we will see, for sufficiently large values of $k$ the minimal speed
is identical to the value $c_\infty$.
\subsection{Analysis of the decay-rates at $+\infty$}
\begin{figure}[h]
\centering
\includegraphics[height=75mm]{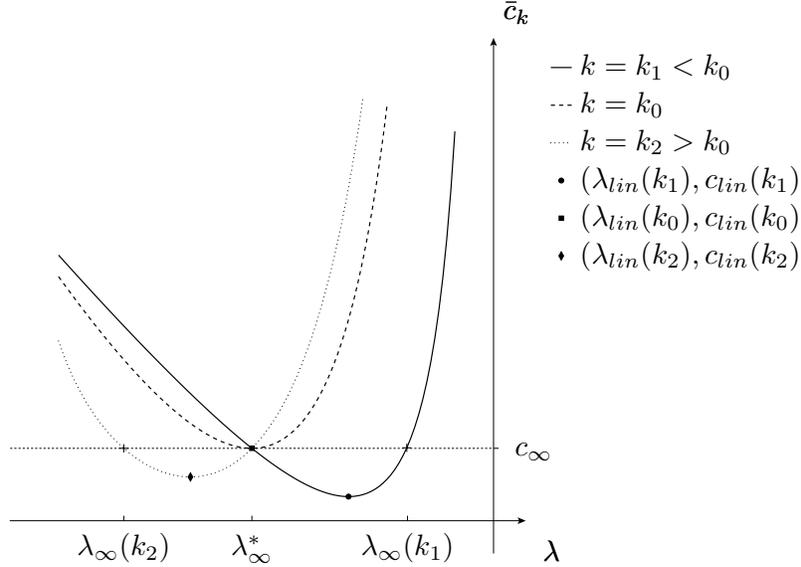}
\caption{Graph of the functions $\bar{c}_k$ for different values of $k$
  and the decay rates $\lambda_\infty(k)$ of the travelling waves with
  speed $c_\infty$ as found in Proposition \ref{prop-red}.}
\label{fig:dec}
\end{figure}

In this section we further investigate the deacy of travelling
waves to zero as $x$ tends to $+\infty$. With this aim we analyze the
functions $\bar{c}_k$ defined in Lemma \ref{lem:curve-lambda-c}:
$\bar{c}_k(\lambda)$ is the
speed of a travelling wave with decay rate $\lambda$ at $+\infty$.
\begin{remark}\label{rem-decay-red}
Corresponding to the reduction \eqref{eq-red-u}, \eqref{eq-red-w} of the
travelling-wave system \eqref{eq-tw-u}, \eqref{eq-tw-w}, we find that
the equation \eqref{eq-lambda} for the possible decay rates factorises
for the speed $c_\infty=(\gamma(1+\gamma))^{-1/2}$. The value
\begin{gather}
  \lambda_\infty^* \,:=\, -\frac{1}{\gamma c_\infty}\,=\,
  -\sqrt{\frac{1+\gamma}{\gamma}} \label{eq:def-lambda-infty-*}
\end{gather}
is for all $k>0$ a negative root of \eqref{eq-lambda} with $c=c_\infty$.
The decay rate of $u_0,w_0$ is given by the other negative root of this equation, which is
the value $\lambda_\infty(k)$ defined in \eqref{eq:dec-u-0}.
\end{remark}
\begin{lemma}\label{lem:decay}
For $k>0$ consider the values $\lambda_{\lin}(k), c_{\lin}(k)$ as
defined in Lemma \ref{lem:curve-lambda-c} and the values
$\lambda_{\infty}(k),\lambda_\infty^*$ as given in \eqref{eq:dec-u-0},
\eqref{eq:def-lambda-infty-*}. Then there exists a unique $k_0>0$, which
is is given explicitely by $k_0=(1+2\gamma)\gamma^{-2}$, such that
\begin{gather}
  \lambda_{\lin}(k_0)\,=\, \lambda_\infty^*.\label{eq:k0}
\end{gather}
Moreover, for $0<k_1<k_0<k_2$
%\begin{gather}
%  c_{\lin}(k_0)\,=\,c_\infty\label{eq:k0-c}
%\end{gather}
\begin{alignat}{3}
  &\lambda_\infty(k_1) \,&&>\,\lambda_{\lin}(k_1)\,&&>\,\lambda_\infty^*
  ,\label{eq:ineqs-k1}\\
  &\lambda_\infty(k_0)
  &&=\lambda_{\lin}(k_0)\,&&=\,\lambda_\infty^*,\label{eq:ineqs-k0}\\ 
  &\lambda_\infty(k_2) \,&&<\,\lambda_{\lin}(k_2)\,&&<\,\lambda_\infty^*
  \label{eq:ineqs-k2}
\end{alignat}
hold, see Figure \ref{fig:dec}.
\end{lemma}
\begin{proof}
We have proved in Lemma \ref{lem:curve-lambda-c} that the functions
$\bar{c}_k$ attain their minimum $c_{\lin}(k)$ at a unique 
value $\lambda_{\lin}(k)$,
\begin{align}
  &c_{\lin}(k)\,=\, \bar{c}_k(\lambda_{\lin}(k)),\label{eq:ck-min}\\
  &\bar{c}_k(\lambda)\,>\,c_{\lin}(k)\quad\text{ for all }\lambda<0,
  \lambda\,\neq\,\lambda_{\lin}(k).\label{eq:ck-min-uni}
\end{align}
For convenience we recall that
\begin{align}
  3\lambda_{\lin}^2 \,&=\, -(1+\gamma k)+2\sqrt{(1+\gamma
  k)^2+3k},\quad \lambda_{\lin}<0,\label{eq:recap-lambda-lin}\\ 
  \lambda_{\lin} c_{\lin} \,&=\, (1+\gamma k) -\sqrt{(1+\gamma
  k)^2+3k}.\label{eq:recap-c-lin} 
\end{align}
By Proposition \ref{prop-red}, Remark \ref{rem-decay-red} and the
definition of $\bar{c}_k$ in Lemma \ref{lem:curve-lambda-c}
\begin{gather}
  \bar{c}_k(\lambda_\infty^*)\,=\,\bar{c}_k(\lambda_\infty(k)) \,=\,
  c_\infty\quad\text{ for all }k>0 \label{eq:eqs-ck}
\end{gather}
holds.
Next we see from
\eqref{eq:recap-lambda-lin} that $\lambda_{\lin}(k)$ is strictly
decreasing in $k$ 
and that
\begin{align*}
  \lim_{k\to 0} \lambda_{\lin}(k)\,&=\, -\frac{1}{3}\sqrt{3}\,>\, -1,\qquad
  \lim_{k\to\infty} \lambda_{\lin}(k) \,=\, -\infty.
\end{align*}
Since $\lambda_\infty^*<-1$ by \eqref{eq:def-lambda-infty-*} there is a
unique value $k_0$ such that \eqref{eq:k0} holds.
By \eqref{eq:k0}, \eqref{eq:ck-min} and \eqref{eq:eqs-ck} we deduce that
\begin{gather}
  c_{\lin}(k_0)\,=\, \bar{c}_{k_0}(\lambda_{\lin}(k_0))
  \,=\, \bar{c}_{k_0}(\lambda^*_\infty)\,=\,
  c_\infty,\label{eq:k0-c}
\end{gather}
and by \eqref{eq:eqs-ck}, \eqref{eq:k0-c} that
\begin{gather}
  c_{\lin}(k_0)\,=\,\bar{c}_{k_0}(\lambda_\infty(k_0)).
  \label{eq:k0-lambda}
\end{gather}
By \eqref{eq:ck-min}, \eqref{eq:ck-min-uni} and \eqref{eq:k0} this yields
\begin{gather*}
  \lambda_\infty(k_0)\,=\, \lambda_\infty^*.
\end{gather*}
Finally one derives from \eqref{eq:k0}, \eqref{eq:recap-c-lin} and
\eqref{eq:k0-c} that $k_0$ is given by
\begin{gather*}
  k_0\,=\, \frac{1+2\gamma}{\gamma^2}.
\end{gather*}
To prove the inequalities \eqref{eq:ineqs-k1}, \eqref{eq:ineqs-k2} we
first observe that \eqref{eq:eqs-ck} implies that
for $k\neq k_0$  there is a $\lambda$ between $\lambda_{\infty}^*$ and
$\lambda_{\infty}(k)$ such that $\bar{c}_k^\prime(\lambda)=0$. By
\eqref{eq:mon-ck-1}, \eqref{eq:mon-ck-2} we conclude that
$\lambda=\lambda_{\lin}(k)$ and deduce that $\lambda_{\lin}(k)$
lies between  $\lambda_{\infty}^*$ and
$\lambda_{\infty}(k)$. Since $\lambda_{\lin}(k)$ is monotonically
decreasing in $k$ we obtain from \eqref{eq:k0} that
\begin{gather*}
  \lambda_\infty(k_1) \,>\,\lambda_\infty^*,\qquad
  \lambda_\infty(k_2) \,<\,\lambda_\infty^*,
\end{gather*}
which proves
\eqref{eq:ineqs-k1}, \eqref{eq:ineqs-k2}.
\end{proof}
The conclusions of Lemma \ref{lem:decay} are illustrated in Figure
\ref{fig:dec}. 
\\

By the previous results we can now compare the decay of two different
travelling waves as $x\to +\infty$, similarly as in Corollary
\ref{cor:decay-} for the convergence to unity as $x\to -\infty$. 
\begin{lemma}\label{lem:comp-decay}
Let $k>k_0$ and assume that $(c,u,w)$ is a monotone travelling wave with speed
$c<c_\infty$. Then, as $x\to +\infty$,
$u(x),w(x)$ decay slower to zero
than does the travelling wave $u_0(x),w_0(x)$ obtained in Proposition
\ref{prop-red}. 
\end{lemma}
\begin{proof}
By \eqref{eq:mon-ck-1} the function $\bar{c}_k$ is monotonically
decreasing for $\lambda<\lambda_{\lin}(k)$. From \eqref{eq:ineqs-k2} and
\eqref{eq:eqs-ck} we therefore deduce that
\begin{gather*}
  \bar{c}_k(\lambda)\,>\, \bar{c}_k(\lambda_\infty(k))\,=\,c_\infty\quad\text{ for all }
  \lambda<\lambda_\infty(k).
\end{gather*}
Since we have assumed that $c<c_\infty$ this implies that the travelling
wave $(c,u,w)$ decays with a rate $\lambda>\lambda_\infty(k)$. On the
other hand, by Proposition \ref{prop-red}, the rate of the exponential
deacy of $u_0,w_0$ as $x\to +\infty$ is given by
$\lambda_\infty(k)$.
\end{proof}
\subsection{The nonlinear selection regime}
In this section we prove that for $k>k_0$, where $k_0$ is given in Lemma
\ref{lem:decay}, the minimal speed of 
travelling waves for the reaction-diffusion system 
\eqref{rd-u}, \eqref{rd-w} is identical to the minimal speed of travelling waves for the
Stefan-like limit problem \eqref{st-w=1}-\eqref{st-Gamma}. In
particular there is nonlinear selection of the minimal speed for $k>k_0$.
This result follows from a comparison principle which is formulated in the 
next Theorem. In general, invariant region arguments do not apply for elliptic systems,
but here a shift parameter is chosen to play the role of
the time parameter in the proof of comparison principles for parabolic systems.
%A similar trick is used in the `moving plane technique', see \cite{Ber}. 
\begin{theorem}\label{the-sel}
Let $(c_1,u_1,w_1)$ and $(c_2,u_2,w_2)$ be two monotone travelling waves and assume that $c_1<c_2$.
Let $\lambda_1,\lambda_2$ denote the decay rates at $+\infty$ of $u_1,w_1$ and $u_2,w_2$
respectively. Then $\lambda_1\leq\lambda_2$ holds and, as $x$ tends to infinity, $u_1(x),w_1(x)$
cannot converge exponentially slower to zero than $u_2(x),w_2(x)$ do.

In particular, a travelling wave has minimal speed if and only if its
decay rate at $+\infty$ is the minimal one among all travelling waves.
\end{theorem}
\begin{proof}
Let us assume that $\lambda_1>\lambda_2$ holds.
Since $c_1<c_2$ the travelling wave $(u_1,w_1)$ converges by Corollary \ref{cor:decay-} 
faster to $(1,1)$ as \mbox{$x\to -\infty$} than does $(u_2,w_2)$. In particular
\begin{eqnarray*}
  (u_1,w_1)(x) &>& (u_2,w_2)(x) \quad\text{ for sufficiently small }x<0.
\end{eqnarray*}
Since we have assumed that $\lambda_1>\lambda_2$, the decay of
$(u_1,w_1)$ at $+\infty$ is slower than the decay of $(u_2,w_2)$ and we
deduce from Lemma \ref{lem-decay-infty+}
that
\begin{eqnarray*}
  (u_1,w_1)(x) &>& (u_2,w_2)(x)\quad\text{ for sufficiently large }x>0.
\end{eqnarray*}
This implies that there is a shift $x_0$, such that for $(\tilde{u}_1,\tilde{w}_1)=(u_1,w_1)(.+x_0)$ 
\begin{gather*}
        \tilde{u}_1 \,\geq\, u_2,\qquad \tilde{w}_1\,\geq\, w_2\quad\text{ on }\R
\end{gather*}
holds and such that there exists a $x_1\in\R$ with
\begin{gather*}
        \tilde{u}_1(x_1)\,=\,u_2(x_1)\quad \text{ or }\quad\tilde{w}_1(x_1)\,=\,w_2(x_1).
\end{gather*}
From equations \eqref{eq-tw-u}, \eqref{eq-tw-w} we obtain for $U:=u_2-\tilde{u}_1, W:=w_2-\tilde{w}_1$
\begin{eqnarray}\label{eq-c,c0-U}
        0 &=& U^{\prime\prime}+cU^\prime - (c_1-c_2)u_2^\prime-U+W -\gamma k (1-w_2)U +\gamma k \tilde{u}_1W,\\
        0 &=& cW^\prime -(c_1-c_2)w_2^\prime +k(1-w_2)U-k\tilde{u}_1W.\label{eq-c,c0-W}
\end{eqnarray}
Assume that $\tilde{u}_1(x_1)=u_2(x_1)$, which gives
\begin{gather*}
        U(x_1) \,=\, \max_x U(x)\,=\,0
\end{gather*}
and $U^\prime(x_1)=0, U^{\prime\prime}(x_1)\leq 0$. Then \eqref{eq-c,c0-U} yields
\begin{gather*}
        0\,\leq\, -(c_1-c_2)u_2^\prime(x_1)\,<\,0,
\end{gather*}
which is a contradiction.\\
If $\tilde{w}_1(x_1)=w_2(x_1)$ then
\begin{gather*}
        W(x_2) \,=\, \max_x W(x)\,=\,0
\end{gather*}
implies $W^\prime(x_2)=0$ and \eqref{eq-c,c0-W} yields
\begin{gather*}
        0\,\leq\, -(c_1-c_2)w_1^\prime(x_2)\,<\,0,
\end{gather*}
which is also a contradiction. Thus we deduce that
$\lambda_1\leq\lambda_2$.

The final conclusion of the Theorem follows now  by a
contradiction argument.
\end{proof}
The comparison principle Theorem \ref{the-sel} and Lemma \ref{lem:comp-decay} imply
that for $k>k_0$ no monotone travelling wave  exists with lower speed than $c_\infty$.
\begin{corollary}\label{cor:minspeed}
For $k\geq k_0$ the minimal speed of travelling waves is given by
\begin{eqnarray*}
  c_{\min}(k) &=& c_\infty.
\end{eqnarray*}
\end{corollary}
\begin{proof}
Assume there is a monotone travelling wave $(c,u,w)$ with $c<c_\infty$ and let $(c_\infty,u_0,w_0)$ be the monotone travelling wave
with speed $c_\infty$ which we have found in Proposition \ref{prop-red}.
By Lemma \ref{lem:comp-decay} the functions $(u,w)$ decay  slower to zero at $+\infty$
than $u_\infty,w_\infty$ do, which is a contradiction to Theorem \ref{the-sel}.
Therefore $c_{\min}(k)\geq c_\infty$ holds and, recalling \eqref{eq-upper-bound}, the conclusion
follows.
\end{proof}
%%===================================================
% convergence of travelling waves
%===================================================
\section{Travelling waves in the fast-degradation-rate limit}
\label{sec-conv}
We complete our investigations by proving that the travelling waves of \eqref{rd-u},
\eqref{rd-w} are, for large values of $k$, close to a travelling waves
of the limit problem \eqref{st-w0}-\eqref{st-Gamma}.
\begin{proposition}\label{prop:speed-lim-k}
Let $(c,u_k,w_k)$, $k\in\N$, be a sequence of monotone travelling waves
for \eqref{rd-u}, \eqref{rd-w} with speed $c\geq c_\infty$ and
\begin{gather}
  u_k(0)\,=\, \frac{1}{2}\quad\text{ for all }k>0. \label{eq:fix-u0}
\end{gather}
Then, as $k$ tends to infinity,
\begin{alignat}{2}
  u_k\,&\to\, U\quad&&\text{ uniformly, }\label{eq:conv-uk-U-1}\\
  u_k\,&\to\, U\quad&&\text{ in }C^{0,\beta}_{\loc}(\R),\text{ for all
  }0<\beta<1,\label{eq:conv-uk-U-2}\\ 
  w_k\,&\to\, W\quad&&\text{ in }L^p_{\loc}(\R),\text{ for all }1\leq
  p<\infty,\label{eq:conv-wk-W-1}
\end{alignat}
where $U,W$ is the unique travelling-wave solution of the limit
problem with speed $c$ and $U(0)=1/2$.
\end{proposition}
\begin{proof}
We recall that $u_k^\prime,w_k^\prime<0$, that $0<u_k,w_k<1$ and
that by \eqref{eq:bound-u-eps-prime}, \eqref{eq:bound-w-eps-prime}
\begin{gather}\label{est-tw-1}
        -u_k^\prime \,<\,\mu
\end{gather}
holds uniformly in $k$, where $\mu$ was defined in \eqref{eq:def-mu}.
This yields the existence of a subsequence $k_i\to\infty\ (i\to\infty)$
and monotone decreasing functions $u,w$ with 
$0\leq u,w\leq 1$ and $u\in H^{1,\infty}_{\loc}(\R)$ such that
\begin{alignat}{3}
        u_{k_i} &\to u \quad&&\text{ in } C^{0,\beta}_{\loc}(\R)&&\text{
        for all }0<\beta<1,\label{conv-tw-u}\\ 
        u_{k_i}^\prime &\schwto u^\prime \quad&&\text{ in
        }L^p_{\loc}(\R)&&\text{ for all }1\leq p <\infty, 
        \label{conv-tw-u-prime}\\
        u_{k_i} &\to u \quad&&\text{ uniformly in }\R&&\text{ for all
        }1\leq p <\infty, 
        \label{conv-tw-u-uniform}\\
        w_{k_i} &\to w \quad&&\text{ in }L^p_{\loc}(\R)&&\text{ for all
        }1\leq p <\infty,\label{conv-tw-w},
\end{alignat}
Integrating equation \eqref{eq-tw-w} over $\R$, we obtain
\begin{eqnarray}
        \int_{-\infty}^\infty u_{k_i}(1-w_{k_i}) &=&
        \frac{c}{k_i},\label{bound-L1-kterm} 
\end{eqnarray}
and by Fatou's Lemma we see that
\begin{gather*}
  \int_{-\infty}^\infty u(1-w)\,\leq\,
  \liminf_{i\to\infty}\int_{-\infty}^\infty u_{k_i}(1-w_{k_i})\,=\,0, 
\end{gather*}
which implies that
\begin{eqnarray}
  u(1-w) &=& 0\quad\text{ almost everywhere in }\R.\label{eq:tw-lim-00}
\end{eqnarray}
By \eqref{eq:conv-uk-U-2}, \eqref{conv-tw-u} we obtain that $u(0)=1/2$
and $u(x)\geq 1/2$ for $x<0$. From \eqref{eq:tw-lim-00} 
we deduce that
\begin{eqnarray*}
  w(x) &=& 1 \quad\text{ for }x<0.
\end{eqnarray*}
The equations \eqref{eq-tw-u} and \eqref{est-tw-1},
\eqref{bound-L1-kterm} yield the estimate 
\begin{eqnarray}
        \int_{-\infty}^0 (1-u_{k_i}) &=& -u_{k_i}^\prime(0) -
        cu_{k_i}(0) +c +\int_{-\infty}^0 (1-w_{k_i}) 
        +\int_{-\infty}^0\gamma ku_{k_i}(1-w_{k_i})
        \notag\\
        &\leq& C(c,\gamma)\label{est-tw-1-u}
\end{eqnarray}
and we obtain that $(1-u)\in L^1\big((-\infty,0)\big)$ and $u(x)\to 1$
as $x\to -\infty$. 
By \eqref{eq-tw-u} and \eqref{est-tw-1}, \eqref{bound-L1-kterm}
\begin{gather}
        \int_0^\infty (u_{k_i}-w_{k_i}) \,=\, -u_{k_i}^\prime(0)
        -\frac{c}{2} -\gamma k_i \int_0^\infty u_{k_i}(1-w_{k_i}) 
        \,\geq\, -C(\gamma,c)\label{est-tw-u}
\end{gather}
holds and we further deduce from \eqref{bound-L1-kterm} that
\begin{gather}\label{est-tw-L1-sum}
        \int_0^\infty (u_{k_i} -\frac{1}{2}w_{k_i}) \,\leq\,
        \int_0^\infty u_{k_i}(1-w_{k_i})\,\leq\, \frac{c}{k_i}, 
\end{gather}
which gives, substracting \eqref{est-tw-u},
\begin{eqnarray}\label{est-tw-L1-w}
        \int_0^\infty w_{k_i} &\leq& C(\gamma,c)
\end{eqnarray}
for $k>1$. Further we find from \eqref{est-tw-L1-sum} that
\begin{eqnarray}\label{est-tw-L1-u}
        \int_0^\infty u_{k_i} &\leq& C(\gamma,c)
\end{eqnarray}
holds. By Fatou's Lemma this implies that $u,w\in L^1((0,\infty))$ and
$u(x),w(x)\to 0$ as $x\to\infty$. By \eqref{conv-tw-u} and since the
limits as $x\to\pm\infty$ of $u_{k_i}, U$ coincide, we deduce
\eqref{eq:conv-uk-U-1} from \cite[Lemma 2.4]{odo}.\\
The equations \eqref{eq-tw-u}, \eqref{eq-tw-w} yield that $u_{k_i},
w_{k_i}$ satisfy 
\begin{eqnarray*}
        0 &=& \int_\R \big(u_{k_i}^\prime + c(u_{k_i}+\gamma
        w_{k_i})\big)\eta^\prime + (u_{k_i}-w_{k_i})\eta 
\end{eqnarray*}
for all $\eta\in C^\infty_c(\R)$ and, according to
\eqref{conv-tw-u}-\eqref{conv-tw-w}, we can pass in this 
equation to the limit $i\to\infty$. This yields
\begin{eqnarray}
        0 &=& \int_\R \big(u^\prime + c(u+\gamma w)\big)\eta^\prime +
        (u-w)\eta.\label{eq:tw-lim-1} 
\end{eqnarray}
Since $u,w$ are monotone decreasing from unity to zero and satisfy
\eqref{eq:tw-lim-00}, we deduce that there is a $x_0\in\R$ such that
\begin{alignat*}{2}
  u(x)\,&=\,0&\qquad&\text{ for }x\geq x_0,\\
  u(x)\,&>\,0,\, w(x)\,=\,1&\quad&\text{ for }x<x_0.
\end{alignat*}
Therefore \eqref{eq:tw-lim-1} yields that
\begin{alignat*}{2}
  0\,&=\, u^{\prime\prime}+cu^\prime -u +1&\quad&\text{
    for }x<x_0,\\
  0\,&=\, \gamma c w^\prime + w&\quad&\text{ for }x>x_0,
\end{alignat*}
and that the jump condition
\begin{align*}
  \gamma c (1-\lim_{x\downarrow x_0}w(x)) &\,=\,
    -\lim_{x\uparrow x_0}u^\prime(x)
\end{align*}
has to be satisfied. This shows that $u,w$ is a travelling wave with
speed $c$ of the limit system \eqref{st-w=1}, \eqref{st-u=0}.
\end{proof}
%=======================================================
% conclusions
%=======================================================
\section{Conclusions}\label{sec-concl}
We conclude our investigations with a brief summary and discussion of
our results. 
\subsection{Slow and fast decay at the threshold $c=c_\infty$}
\begin{figure}[h]
\centering
\includegraphics[height=75mm]{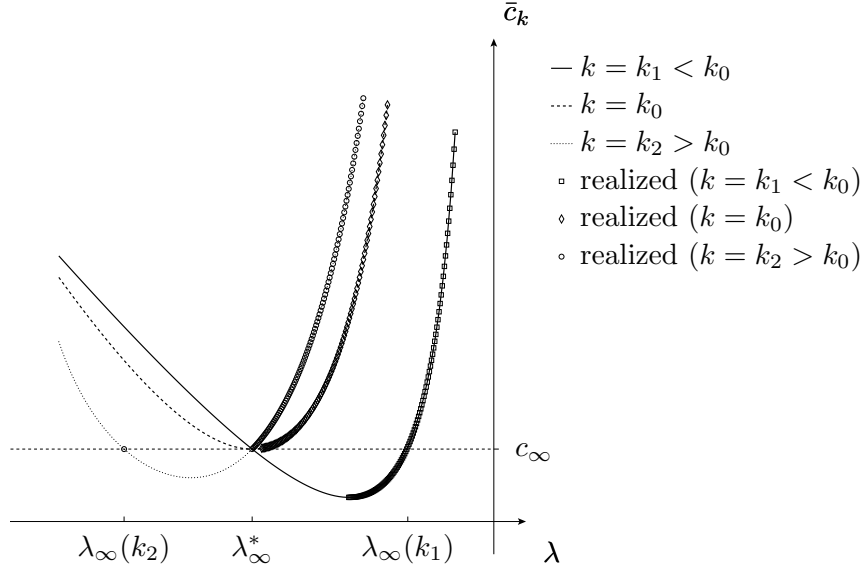}
\caption{Graph of the functions $\bar{c}_k$ for different values of $k$ and the
  points on the graph that are realized by a travelling wave.}
\label{fig:dec2}
\end{figure}

Our results on the existence of travelling waves for the system
\eqref{rd-u}, \eqref{rd-w} and the selection mechanism of the minimal
speed are summarized in Figure \ref{fig:dec2}. As in the preceding
figure we have plotted the functions $\bar{c}_k$ which give for a
travelling wave with speed $c$ the possible rates $\lambda$ of the
exponential decay to zero at $x\to +\infty$. By the circles, squares and
diamonds in
Figure \ref{fig:dec2} we have indicated the decay rates which in fact
are realized by a travelling wave: presuming that there is linear
selection for $k<k_0$ travelling waves exists
for all speeds $c\geq c_{\lin}$.
By Theorem \ref{the-sel} the decay rates of these
travelling waves corresponds to values on the increasing branch of
$(\lambda,c_k(\lambda))$.
This behaviour changes for $k> k_0$: the decay rate of $u_0,w_0$ is
on the decreasing branch of the solution curve. By Corollary \ref{cor:minspeed}
$c_\infty$ is the minimal speed and the decay rates of travelling waves with larger speeds
are on the increasing branch.
\\

To give an explanation of what happens in the nonlinear selection regime
$k>k_0$ if the speed $c$ falls below the minimal speed
$c_{\min}(k)=c_\infty$ we consider the linearization of \eqref{rd-u},
\eqref{rd-w} at $(u^\prime,u,w)=(0,0,0)$, which was given in \eqref{eq:xi-lin}. For all speeds
$c$ in a neighbourhood of $c_\infty$ the linearized system has two
negative eigenvalues; the stable manifold for \eqref{rd-u},
\eqref{rd-w} at $(u^\prime,u,w)=(0,0,0)$ is two-dimensional. One checks that the
eigenspace corresponding to a negative eigenvalue intersects with the
set $\{z=(z_1,z_2,z_3) : z_1<0, z_2>0, z_3>0 \}$. A monotone travelling
wave exists for $c>c_0$ and approaches $(0,0,0)$ tangentially to the
eigenspace corresponding to the larger negative eigenvalue (`slow
decay'), see Figure \ref{fig:dec2}. We expect also for $c<c_\infty$ an orbit
$(u^\prime,u,w)$ connecting $(0,1,1)$ with $(0,0,0)$. Such an orbit will
also converge to $(0,0,0)$ tangentially to the eigenspace corresponding
to the larger negative eigenvalue but comes from the `wrong' side, 
taking negative values for $u,w$. For the threshold $c=c_\infty$ there
still exists a monotone travelling wave; this travelling wave
approaches $(0,0,0)$ tangentially to the eigenspace of the smaller
negative eigenvalue (`fast decay'). 

\subsection{Discussion}
Travelling waves often determine the long-time behaviour of solutions for
the initial-value problem for \eqref{rd-u}, \eqref{rd-w}. Typically,
for solutions with sufficiently fast decaying initial data, the
propagation speed of pertubations from the unstable equilibrium 
is given by the minimal speed of travelling waves. The proof of
such a result, as well as the uniqueness of travelling waves,
for the system \eqref{rd-u}, \eqref{rd-w} is not in the
scope of the present article. Nevertheless our result that the
travelling wave with 
minimal speed has the fastest decay supports that conjecture.

For the reaction-diffusion system \eqref{rd-u}, \eqref{rd-w} in
arbitrary space-dimension $n$ our analysis yields a family of
supersolutions: Consider for  $c\geq
c_\infty$ a travelling 
wave $u,w$ with speed $c$ and an arbitrary real number $r\in \R$. The
functions $\bar{u},\bar{w}$ defined by 
\begin{gather*}
  \bar{u}(x,t)\,:=\, u(|x|-ct-r),\quad \bar{w}(x,t)\,:=\, w(|x|-ct-r)
\end{gather*}
satisfy
\begin{align}
  \partial_t {u} - \Delta{u} + {u} - {w} +\gamma k {u}(1-{w})\,&=\,
  -\frac{n-1}{|x|} \bar{u}^\prime,\label{eq:n-dim-rd}\\
  \partial_t {w} - k {u}(1-{w})\,&=\, 0.\notag
\end{align}
Since $\bar{u}^\prime$ is negative $(\bar{u},\bar{w})$ is
a supersolution.
In order to construct a subsolution one has to control the
 dimension-depending correction term in \eqref{eq:n-dim-rd}.

In view of the applications, the robustness of the wave-speed to changes
in the parameter values is a valuable feature: see the explicit formula
\eqref{eq-c-infty} for the minimal speed and Proposition
\ref{prop:speed-lim-k}. For the mathematical
analysis of reaction-diffusion  \emph{systems} and the selection of the
minimal speed the model that we have derived is a good paradigm. We
prove that nonlinear selection occurs and determine explicitly the
minimal speed. One crucial ingredient
is the exact first integral obtained in Lemma
\ref{lem:red}, which is a special property of the system \eqref{rd-u},
\eqref{rd-w}. Other results, in particular the comparison principle
Theorem \ref{the-sel} and the observation that the fastest decay is
realized by a travelling wave with minimal speed, can be extended to general
monotone systems. 
%%%%%%%%%%%%%%%%%%%%%%%%
% bibliography
%%%%%%%%%%%%%%%%%%%%%%%%
\newcommand{\etalchar}[1]{$^{#1}$}

\end{document}